\documentstyle[graphicx,amsfonts,amssymb,epstopdf,11pt]{article}

\newcommand{\dsp}{\displaystyle}

\newcommand{\sgn}{\mbox{\rm \small  sgn}}

\newtheorem{theorem}{Theorem}[section]
\newtheorem{lemma}[theorem]{Lemma} 
 
\newtheorem{corollary}[theorem]{Corollary} 
\newtheorem{remark}[theorem]{Remark} 
\newtheorem{proposition}[theorem]{Proposition}

\title{\bf Global in Time Existence of  Self-Interacting Scalar Field in  De~Sitter Spacetimes }
\author{\bf Anahit Galstian, Karen Yagdjian\footnote{E-mail addresses: anahit.galstyan@utrgv.edu, karen.yagdjian@utrgv.edu}\\
{}\\
{\small School of Mathematical and Statistical Sciences,   University of Texas Rio Grande Valley},\\
{\small 1201 West University Dr.,
 Edinburg, TX 78539, U.S.A.}}

\begin{document}

\date{}
\maketitle
\thispagestyle{empty} 
\vspace{-0.3cm}

\medskip

\addtocounter{section}{-1}
\renewcommand{\theequation}{\thesection.\arabic{equation}}
\setcounter{equation}{0}
\pagenumbering{arabic}
\setcounter{page}{1}
\thispagestyle{empty}

\begin{abstract}
We prove the  existence of a global in time    solution of the semilinear Klein-Gordon equation in the de~Sitter space-time.
 
\end{abstract}

\section{Introduction}
\label{S0}

In this paper we prove the  existence of a global in time    solution of the semilinear Klein-Gordon equation in the de~Sitter space-time. 
\smallskip

The metric $g$ in the de~Sitter   
space-time   is defined as follows, 
$g_{00}=  g^{00}= -  1  $, $g_{0j}= g^{0j}= 0$, $g_{ij}(x,t)=e^{2t}    \sigma  _{ij} (x) $, 
$i,j=1,2,\ldots,n$, where $\sum_{j=1}^n\sigma  ^{ij} (x) \sigma _{jk} (x)=\delta _{ik} $,   and $\delta _{ij} $ is  Kronecker's delta. The metric $\sigma  ^{ij} (x) $ describes the time slices. 
In the quantum field theory    the matter fields are described by a function  $\psi $  that must satisfy  equations of motion.
In the case of a massive scalar field, the equation of motion is   the  semilinear Klein-Gordon equation generated by the metric $g$:
\begin{eqnarray*}
\square_g \psi = m^2 \psi   +   V'_\psi (x,\psi  ) \,.
\end{eqnarray*}
Here $m$ is a physical mass of the particle. In physical
terms this equation describes a local self-interaction for a scalar particle.  A typical example of a potential function  would be  $V(\psi  )=\psi  ^4$.
The semilinear equations are also commonly used models for general nonlinear problems. 
\smallskip

In the present paper we are going to extend  the  small data global existence result of \cite{JMAA2012}  for the spatially  flat de~Sitter space-time to the    de~Sitter space-time
with the time slices being Riemannian manifolds. 
\smallskip

To formulate the main theorem of this paper we need the following description of the nonlinear term. Let $B^{s,q}_{p}$ be the Besov space. \\

\noindent {\bf Condition ($\mathcal L$).} {\it The function $F$ is
said to be Lipschitz continuous with exponent $\alpha \geq 0 $ in the space $B^{s,q}_{p}$   if there is a constant \,
$C \geq 0$ \, such that
\begin{equation}
\label{calM}
 \hspace{-0.4cm} \|  F(x,\psi  _1 (x))- F(x,\psi  _2(x) ) \|_{B^{s,q}_{p}}  \leq C
\| \psi  _1 -  \psi _2   \|_{B^{s,q}_{p'}}
\Big( \|  \psi _1  \|^{\alpha} _{B^{s,q}_{p'}}
+ \|  \psi _2   \|^{\alpha} _{B^{s,q}_{p'}} \Big)
\end{equation}
for all \,\,  $\psi _1,\psi  _2 \in  B^{s,q}_{p'}$, where $1/p+1/p'=1 $. }

For the   case of $B^{s,2}_{2}= H_{(s)} ({\mathbb R}^n)$, define the complete metric space
\[
X({R,s,\gamma})  := \{ \psi  \in C([0,\infty) ; H_{(s)} ({\mathbb R}^n)) \; | \;
 \parallel  \psi   \parallel _X := \sup_{t \in [0,\infty) } e^{\gamma t}  \parallel
\psi  (x ,t) \parallel _{H_{(s)} ({\mathbb R}^n)}
\le R \}\,,
\]
$\gamma  \geq 0$, with the metric
\[
d(\psi _1,\psi _2) := \sup_{t \in [0,\infty) }  e^{\gamma t}  \parallel  \psi _1 (x , t) - \psi _2 (x ,t) \parallel _{H_{(s)} ({\mathbb R}^n)}\,.
\]
We denote  ${\mathcal B}^\infty $ the space of all $C^\infty ({\mathbb R}^n)$ functions with uniformly bounded derivatives of all orders. 
\begin{theorem}
\label{T0.1}
Let  \, $ A(x,\partial_x) = \sum_{|\alpha |\leq 2} a_{\alpha }(x) \partial_x^\alpha $ \, be a second order negative elliptic differential operator with real coefficients \, $a_{\alpha } \in {\mathcal B}^\infty$.   Assume that the nonlinear term $F(u)$  is  Lipschitz continuous with exponent $\alpha > 0 $ in the  space $H_{(s)} ({\mathbb R}^n)$, $ s > n/2\geq 1$, and $F(0)=0$. 
Assume also that $m \in (0, \sqrt{n^2-1}/2]\cup [n/2,\infty) $.  
Then  there exists $\varepsilon _0>0 $ such that, for every given functions $\psi _0 ,\psi _1 \in H_{(s)} ({\mathbb R}^n) $,  such that 
\[
 \| \psi _0   \|_{ { H}_{(s)} ({\mathbb R}^n)}
+  \|\psi _1  \|_{ { H}_{(s)} ({\mathbb R}^n)} \leq \varepsilon, \qquad \varepsilon  < \varepsilon_0\,, 
\] 
there exists a global solution $\psi  \in C^1([0,\infty);H_{(s)} ({\mathbb R}^n))$ of the Cauchy problem  
\begin{eqnarray}
\label{NWE}  
&  &
\psi _{tt} +   n  \psi  _t - e^{-2 t} A(x,\partial_x) \psi   +  m^2  \psi  =   F(x, \psi  )\,,\\
\label{ICPHI} \label{0.5}
&  &
\psi (x,0)=\psi _0 (x)\, , \quad \psi _t(x,0)=\psi _1 (x) \,.
 \end{eqnarray} 
The  solution \, $ \psi  (x ,t) $ \, belongs to the space \, $  X({2\varepsilon,s,\gamma})  $, that is,  
\begin{eqnarray*}  
\sup_{t \in [0,\infty)}  e^{\gamma t}  \|\psi  (\cdot ,t) \|_{H_{(s)} ({\mathbb R}^n)}  < 2\varepsilon \,,
\end{eqnarray*}
with  \, $\gamma  $  \, such that  either \,  $0 < \gamma \leq   \frac{1}{\alpha +1 }\left( \frac{n}{2 } -\sqrt{\frac{n^2}{4 }-m^2} \right) $  if \,  $\sqrt{n^2-1}/2\geq  m>0$,  or   we choose $0\leq \gamma _0 <  \frac{n-1}{2}$ if  $m=n/2$ and $0\leq \gamma _0 \leq   \frac{n-1}{2}$ if $m>n/2$,
 then   
$\gamma    \leq  \min \left\{ \gamma _0 \,,\,  \frac{n}{2(\alpha +1)} \right\} $.

If $ m \in ( \sqrt{n^2-1}/2,n/2)$, then for the problem with $\psi _0=0 $ the global solution exists and   belongs 
to $  X({2\varepsilon,s,\gamma}) $, where $\gamma \in (0,  \frac{1}{\alpha +1 }( \frac{n}{2 } -\sqrt{\frac{n^2}{4 }-m^2}) )$.
 \end{theorem}

The range $m \in (\sqrt{n^2-1}/2,n/2) $, which seems to be a forbidden mass interval for the problem with general initial data, can be allowed  if we  change the setting of the problem. Indeed, if we consider the initial value problem with vanishing Cauchy data and with the source term $f$,    then we have the following result for all $ m>0$.
\begin{theorem}
\label{T0.2}
Let  \, $ A(x,\partial_x) = \sum_{|\alpha |\leq 2} a_{\alpha }(x) \partial_x^\alpha $ \, be a second order negative elliptic differential operator with real coefficients \, $a_{\alpha } \in {\mathcal B}^\infty$.   Assume that the nonlinear term $F(u)$  is a  Lipschitz continuous with exponent $\alpha > 0 $ in the  space $H_{(s)} ({\mathbb R}^n)$, $ s > n/2\geq 1$, and $F(0)=0$. 
Assume also that $m >0 $.  
Then  there exists $\varepsilon _0>0 $ such that, for every given function $f \in X({\varepsilon,s,\gamma_{rhs} }) $, 
such that 
\begin{eqnarray*}  
 \sup_{t  \in [0,\infty)} e^{\gamma _{rhs}t}\| f (x ,t) \|_{H_{(s)} ({\mathbb R}^n)}  \leq \varepsilon  <  \varepsilon _0 \,,
\end{eqnarray*}
there exists a global solution $\psi  \in C^1([0,\infty);H_{(s)} ({\mathbb R}^n))$ of the Cauchy problem  
\begin{eqnarray}
\label{NWEf}  
&  &
\psi _{tt} +   n  \psi  _t - e^{-2 t} A(x,\partial_x) \psi   +  m^2  \psi  = f+  F(x,\psi  )\,,\\
\label{ICPHIf} 
&  &
\psi (x,0)=0\, , \quad \psi _t(x,0)=0 \,.
 \end{eqnarray} 
The  solution \, $ \psi  (x ,t) $ \, belongs to the space \, $  X({2\varepsilon,s,\gamma})  $, that is,  
\begin{eqnarray*}  
\sup_{t \in [0,\infty)}  e^{\gamma t}  \|\psi  (\cdot ,t) \|_{H_{(s)} ({\mathbb R}^n)}  < 2\varepsilon \,,
\end{eqnarray*}
with  \, $\gamma  $  \, such that 
\begin{eqnarray*}  
\cases{ 
\dsp \gamma < \frac{1}{\alpha +1 }  \gamma _{rhs}   \hspace{1.5cm}   \quad if \quad m < \frac{n}{2} \quad and     \quad  \gamma _{rhs}\leq \frac{n}{2}-\sqrt{\frac{n^2}{4}-m^2}, \cr
\dsp \gamma < \frac{1}{\alpha +1 } \left( \frac{n}{2}-\sqrt{\frac{n^2}{4}-m^2} \right)   \quad if  \quad m < \frac{n}{2} \quad and     \quad  \gamma _{rhs}> \frac{n}{2}-\sqrt{\frac{n^2}{4}-m^2}, \cr
\dsp \gamma \leq  \min \left\{\gamma _{rhs},  \frac{n}{2(\alpha +1)}\right\}\quad if  \quad m\geq   \frac{n}{2} \quad and     \quad \frac{n}{2} > \gamma _{rhs} ,\cr
\dsp \gamma \leq  \min\left\{ \gamma _{0}, \frac{n}{2(\alpha +1)}\right\} \quad  where \quad \gamma _{0}< \gamma _{rhs}\quad if  \quad m =  \frac{n}{2} \quad and     \quad \frac{n}{2} = \gamma _{rhs},\cr
\dsp \gamma \leq   \frac{n}{2(\alpha +1)}  \hspace{1.5cm}\quad if   \quad m >   \frac{n}{2} \quad and     \quad \frac{n}{2} < \gamma _{rhs},\cr
\dsp \gamma <   \frac{n}{2(\alpha +1)}  \hspace{1.5cm}\quad if   \quad m =   \frac{n}{2} \quad and     \quad \frac{n}{2} < \gamma _{rhs},\cr
\dsp \gamma \leq   \frac{n}{2(\alpha +1)} \hspace{1.5cm}  \quad if  \quad m >  \frac{n}{2} \quad and     \quad \frac{n}{2} = \gamma _{rhs}. 
}
\end{eqnarray*} 
 \end{theorem}
 \medskip

The main tools to prove Theorems~\ref{T0.1},\ref{T0.2} are the following: 1) the integral transform, which produces   representations of the solutions of the linear equation,  
2)  the decay estimates in the Besov spaces, which generate weighted Stricharz estimates, and 3) the fixed point theorem.
\medskip

The kernel of the integral transform  that will be used in this paper (see Section~\ref{S1}) is the following function
\begin{eqnarray*}
E(x,t;x_0,t_0;M)
& :=  & 
 4 ^{-M}  e^{ M(t_0+t) } \Big((e^{-t _0}+e^{-t})^2 - (x - x_0)^2\Big)^{M -\frac{1}{2}} \\
 &  &
\times F\Big(\frac{1}{2}-M   ,\frac{1}{2}-M  ;1;
\frac{ ( e^{-t_0}-e^{-t })^2 -(x- x_0 )^2 }{( e^{-t_0}+e^{-t })^2 -(x- x_0 )^2 } \Big) , \nonumber
\end{eqnarray*}
where   $(x,t) \in D_+ (x_0,t_0)\cup D_- (x_0,t_0) $, $D_+ (x_0,t_0)$ and $D_- (x_0,t_0) $ are the {\it chronological future} and the {\it chronological past}, respectively,   of the point   $(x_0 ,t_0)\in {\mathbb R}^{n+1}$,   \,    $M \in {\mathbb C}$, and $F\big(a, b;c; \zeta \big) $ is the hypergeometric function. The values of the physical mass $m$, which lead to the values of  $M=-k+\frac{1}{2} $, $k=0,1,2,\ldots$, are called in \cite{JMP2013} {\it the knot points}. 
One of these knot points, $m=\sqrt{ n^2-1}/{2}$, presents the only field that obeys the Huygens' principle \cite{JMP2013}.
For these values of the  curved mass $M$ the functions $F(-k,-k;1;z) $, $k=0,1,2,\ldots\,$, are polynomials. 
\medskip

The Klein-Gordon quantum fields on the de~Sitter manifold  with 
imaginary mass, which   take an infinite
set of discrete values as follows 
\begin{eqnarray}
\label{E-M}
m^2=-k(k+n)\,, \quad k=0,1,2,\ldots,
\end{eqnarray}  
present a family of scalar tachyonic  quantum fields.    The nonexistence of a global in time  solution of the semilinear Klein-Gordon massive tachyonic (quantum
fields)  equation in the de~Sitter space-time is proved in \cite{Yag_DCDS}. The conclusion is that the self-interacting  tachyons  in the de Sitter space-time have finite lifespan. More precisely, consider the   semilinear equation
\begin{eqnarray*}
&  &
\psi _{tt} +   n  \psi  _t - e^{-2 t} \Delta  \psi   -  m^2  \psi  =  c|\psi|^{1+\alpha} \,,  
 \end{eqnarray*} 
which is   commonly used model for general nonlinear problems. Then, according to Theorem~1.1~\cite{Yag_DCDS}, if $c\not= 0$, $\alpha >0 $,
and $m \not= 0$, then for every positive numbers $\varepsilon  $ and $s$ there exist  functions $\psi _0 $, 
$\,  \psi _1  \in C_0^\infty ( {\mathbb R}^n)$ such that $ \|\psi _0 \|_{H_{(s)}( {\mathbb R}^n)} + \|\psi _1 \|_{H_{(s)}( {\mathbb R}^n)} \leq \varepsilon $ but the solution $\psi =\psi (x,t) $ with the  initial values (\ref{ICPHI}) blows up in finite time. 
\medskip

The equation, which is  considered in Theorem~\ref{T0.1} and Theorem~\ref{T0.2},
is  more general than the covariant Klein-Gordon equation. Furthermore, these theorems, after evident modification,
 can be applied to the smooth pseudo Riemannian manifold $({\mathcal V},g)$ of dimension $n+1$ and ${\mathcal V}={\mathbb R}\times {\mathcal S} $ with ${\mathcal S}$ an $n$-dimensional orientable smooth manifold
and $g$ is the de~Sitter metric. 
One important example  of the equation on the smooth pseudo Riemannian manifold that is amenable to Theorem~\ref{T0.1} is 
if ${\mathcal S}$ is a non-Euclidean space of constant negative curvature  and then the equation of the problem (\ref{NWE})\&(\ref{ICPHI})
is a non-Euclidean Klein-Gordon equation.
\medskip

Although recently the equations in the de Sitter and anti-de Sitter space-times became the focus of interest for an
increasing number of authors (see, e.g., \cite{Allen-Folacci,Bachelot,BaskinSE,Brevik-Simonsen,Epstein-Moschella,GY_NA,Hintz-Vasy,Moschella,Nakamura,
Vasy_2010} and
the bibliography therein) which investigate those equations from a wide spectrum of perspectives, 
there are very few papers on the semilinear Klein-Gordon equation in the de~Sitter space-time. Here we mention some of them closely related to our main result. 
Baskin~\cite{BaskinSE} discussed small data global energy class solutions for the
scalar Klein-Gordon equation on asymptotically de Sitter spaces,
which are   compact manifolds with boundary.   More precisely, in
\cite{BaskinSE}  the following Cauchy problem is considered for the
semilinear equation
\begin{eqnarray*}
&  &
\square_g u +m^2 u = f(u),\qquad
u(x,t_0)=\varphi_0 (x) \in H_{(1)}({\mathbb R}^n) , \,\,\, u_t(x,t_0)=\varphi_1 (x) \in L^2 ({\mathbb R}^n)\,,
\end{eqnarray*}
where mass is large, $m^2> n^2/4$, $f$ is a smooth function and
satisfies conditions $ |f(u)| \leq c|u|^{\alpha  +1}$, $ |u| \cdot
|f'(u)| \sim |f(u)|$, $f(u)-f'(u)\cdot u   \leq 0$, $  \int_0^u
f(v)dv \geq 0$, and $\int_0^u f(v)dv \sim |u|^{\alpha  +2} $ for
large $|u|$. It is also assumed that $  \alpha= \frac{4}{n-1}$. In
Theorem~1.3 \cite{BaskinSE} the existence of the global solution for
small energy data is stated. (For more references on the
asymptotically de Sitter spaces, see  the bibliography in
\cite{Baskin},
 \cite{Vasy_2010}.)
\smallskip

The case of the massless field $m=0$, that is, the global existence for the semilinear   wave equation  (\ref{NWE}) in the de~Sitter space-time, is still open. 
On the other hand, de~Sitter space-time has different realizations (see, e.g., \cite{Hawking}), that allow one to reduce the problem   to the case of the manifold with  constant curvature. 
For geometric reasons, one can expect better dispersive properties from a linear equation on the manifold with negative constant curvature, and, consequently,    stronger results for a semilinear equation, than
in the Euclidean setting.
\smallskip

The global existence of the solutions of the equation on  the manifold with the time slices being real hyperbolic
spaces ${\mathbb H}^n$ is investigated by   Metcalfe and Taylor in \cite{Metcalfe-Taylor,Metcalfe-Taylor_2012},  and by  Anker,   Pierfelice,  and  Vallarino in \cite{Anker,Anker_2014}. In particular, in \cite{Metcalfe-Taylor} for the equation $(\partial_t^2 -\Delta )u=au^2 $ 
on ${\mathbb R}\times {\mathbb H}^3 $ existence of the  small data global solution is proved. Thus, it is shown for a range of powers that is broader than that known for Euclidean space from the
so-called Strauss conjecture. 
Note that the operators in   those articles \cite{Metcalfe-Taylor,Metcalfe-Taylor_2012,Anker,Anker_2014} have time-independent coefficients. 
\smallskip

The Cauchy problem for the damped linear wave equations 
with a time-dependent propagation speed and dissipations, $u_{tt}-a(t)^2\Delta u + b(t)u_t=0 $, where $a \in L^1(0,\infty) $, is considered by Ebert and Reissig  \cite{Ebert-Reissig}. The analysis  of results of   \cite{Ebert-Reissig} hopefully  can lead to the global existence in the problem for the wave equation in the de~Sitter space-time and can shed a 
light on the  interval $(\sqrt{n^2-1}/2,n/2) $. 
\smallskip

Nakamura   \cite{Nakamura} considered the Cauchy problem for the semi-linear Klein-Gordon equations  in de~Sitter space-time with $n \leq 4 $ and with flat time slices. 
The nonlinear term  is of power type for $n=3,4$, or of exponential type for $n=1,2$. For the power type semilinear term with \, $ \frac{4}{n} \leq  \alpha \leq \frac{2}{n-2} $ Nakamura \cite{Nakamura} proved  the existence of    global solutions in the energy class.
\smallskip

Galstian and Yagdjian~\cite{GY_NA} proved  the existence of    global solutions in the energy class in the case of $n =3, 4 $ and the nonlinear term  is of power type. They considered the equation in the  Friedmann-Lema$\hat{i} $tre-Robertson-Walker   space-times  (FLRW space-times) 
with the time slices being Riemannian manifolds. The Klein-Gordon equation in the Einstein-de~Sitter and de~Sitter space-times are important particular cases discussed in ~\cite{GY_NA}.
\smallskip

The present paper is organized as follows.  In Section~\ref{S1} we describe 
 the  integral transform and,   generated by that transform, representations (from \cite{MN})
for the solutions of the Cauchy problem for the linear equation. Then, in Sections~\ref{S2}-\ref{S4} we  derive 
  the $ B^{s,q}_{p}-B^{s',q}_{p'}$ estimates for the solutions of that equation with large, small and critical (for the Huygens' principle) mass, respectively.  
  These representations are used in the Subsections~\ref{SS3.2cr} -\ref{SS3.3cr} for the derivation of asymptotic expansions.
  The last section, Section~\ref{S5}, is devoted to the solvability of the associated integral equation and to the proof of Theorem~\ref{T0.1} and 
 Theorem~\ref{T0.2}.

\section{Integral Transform}
\label{S1}

\setcounter{equation}{0}

If we introduce the new unknown function $u=e^{nt/2}\psi  $, then the equation (\ref{NWE}) takes the form of the  Klein-Gordon equation 
\begin{equation}
\label{1.1}
u_{tt}
-   e^{-2t} A(x,\partial _x) u  + M^2 u = e^{nt/2}F(x,e^{-nt/2} u) \,,
\end{equation}
where $M^2=m^2 -\frac{n^2}{4}$ is the square of the so-called curved (or effective) mass.

We introduce the following notations.
First, we define a  {\it  chronological future}
$D_+ (x_0,t_0) $
and  a
{\it  chronological past} $D_- (x_0,t_0) $  of the point   $(x_0 ,t_0)$,  $x_0 \in {\mathbb R}^n$, $t_0 \in {\mathbb R}$,
 as follows:
$
D_\pm (x_0,t_0)
  :=
 \{ (x,t)  \in {\mathbb R}^{n+1}  \, ; \,
|x -x_0 | \leq $ $\pm( e^{-t_0} - e^{-t })
\, \} $.
Then, for $(x_0, t_0) \in {\mathbb R}^n\times {\mathbb R}$, $M \in {\mathbb C}$,  we define the function
\begin{eqnarray}
\label{EM}
E(x,t;x_0,t_0;M)
& :=  & 
 4 ^{-M}  e^{ M(t_0+t) } \Big((e^{-t _0}+e^{-t})^2 - (x - x_0)^2\Big)^{M -\frac{1}{2}} \\
 &  &
\times F\Big(\frac{1}{2}-M   ,\frac{1}{2}-M  ;1;
\frac{ ( e^{-t_0}-e^{-t })^2 -(x- x_0 )^2 }{( e^{-t_0}+e^{-t })^2 -(x- x_0 )^2 } \Big) , \nonumber
\end{eqnarray}
where   $(x,t) \in D_+ (x_0,t_0)\cup D_- (x_0,t_0) $  and $F\big(a, b;c; \zeta \big) $ is the hypergeometric function. (For definition of the hypergeometric function, see, e.g., \cite{B-E}.)
When no ambiguity arises, like in (\ref{EM}), we use the notation $x^2:= |x|^2$ for $x \in {\mathbb R}^n $.
Thus, the function $E  $ depends on $r^2= (x- x_0 )^2$, that is $E(x,t;x_0,t_0;M)= E  (r,t;0,t_0;M) $. According to Theorem~2.12~\cite{MN}, the function $  E  (r,t;0,t_0;M) $ solves  the following one dimensional Klein-Gordon equation in the de~Sitter space-time:
 \begin{eqnarray*}
&  &
E_{tt} (r,t;0,t_0;M) - e^{-2t} E_{rr}(r,t;0,t_0;M)-M^2 E(r,t;0,t_0;M) =0\,.
\end{eqnarray*}

The kernels   $K_0(z,t;M) $   and  $K_1(z,t;M) $ are defined by
\begin{eqnarray}
\label{K0M}
K_0(z,t;M)
&  := &
- \left[ \partial_b E(z,t;0,b;M) \right]_{b=0}\,,\\
\label{K1M}
K_1(z,t;M)
& :=  &
  E(z,t;0,0;M)\,.  
\end{eqnarray}

The equation (\ref{1.1}) is said to be an equation with imaginary (real) mass if $M^2 > 0 $ ($-M^2 \geq 0 $); here $M \in {\mathbb C}$.
Assume that $A(x,\partial_x) = \sum_{|\alpha| \leq m} a_\alpha (x)\partial^\alpha_x$, where $ a_\alpha (x) \in C^\infty (\Omega ) $ and $\Omega \subseteq {\mathbb R}^n $. For the  Klein-Gordon equation (\ref{1.1})    we invoke to the following result.

\begin{theorem} \mbox{\rm \cite{MN}}
\label{T1.2}
For $ f \in C (\Omega\times I  )$,\, $ I=[0,T]$, $0< T \leq \infty$, and \, $ \varphi_0 $,  $ \varphi_1 \in C (\Omega ) $,
let  the function\,
$v_f(x,t;b) \in C_{x,t,b}^{m,2,0}(\Omega \times [0,1-e^{-T}]\times I)$\,
be a solution to the   problem
\begin{equation}
\label{1.22}
\cases{ 
 v_{tt} -   A(x,\partial_x) v  =  0 \,, \quad x \in \Omega \,,\quad t \in [0,1-e^{-T}]\,,\cr
v(x,0;b)=f(x,b)\,, \quad v_t(x,0;b)= 0\,, \quad b \in I,\quad x \in \Omega\,,
} 
\end{equation}
and the function \, $  v_\varphi(x, t) \in C_{x,t}^{m,2}(\Omega \times [0,1-e^{-T}])$ \, be a   solution   of the   problem
\begin{equation}
\label{1.23}
\label{1.6}
\cases{
 v_{tt}-   A(x,\partial_x) v =0, \quad x \in \Omega \,,\quad t \in [0,1-e^{-T}]\,, \cr
 v(x,0)= \varphi (x), \quad v_t(x,0)=0\,,\quad x \in \Omega\,.
}
\end{equation}

Then the function  $u= u(x,t)$   defined by
\begin{eqnarray*}
u(x,t)
&  =  &
2   \int_{ 0}^{t} db
  \int_{ 0}^{\phi (t)- \phi (b)}  v_f(x,r ;b) E(r,t;0,b;M)  \, dr  
+ e ^{\frac{t}{2}} v_{\varphi_0}  (x, \phi (t))\\
&  &
+ \, 2\int_{ 0}^{\phi (t)} v_{\varphi_0}  (x, s) K_0( s,t;M)   ds  \nonumber 
+\, 2\int_{0}^{\phi (t) }  v_{\varphi _1 } (x,  s)
  K_1( s,t;M)   ds
, \quad x \in \Omega  , \,\, t \in I\,,
\end{eqnarray*}
with $\phi (t):= 1-e^{-t} $,
 solves the problem
\begin{equation}
\label{1.20}
\cases{
u_{tt} - e^{-2t}A(x,\partial_x) u - M^2 u= f, \quad  x \in \Omega \,,\,\, t \in I,\cr
  u(x,0)= \varphi_0 (x)\, , \quad u_t(x,0)=\varphi_1 (x),\quad x \in \Omega\,.
}
\end{equation}
Here the kernels  $E$, $K_0$ and $K_1$ have been defined in (\ref{EM}), (\ref{K0M}) and (\ref{K1M}), respectively.
 \end{theorem}

We note    that the operator $A(x,\partial_x)$  is of arbitrary order, that is, the equation of (\ref{1.20}) can be an evolution equation, not necessarily hyperbolic.
Then,   the problems  (\ref{1.22}) and (\ref{1.20}) can be   mixed initial-boundary value problems involving the boundary condition.   
Next, we stress that interval $[0,1-e^{-T}]\subseteq [0,1]$, which appears in (\ref{1.22}), reflects the fact that de~Sitter model possesses the horizon \cite{Hawking}; existence of the horizon
in the de~Sitter model is widely used to define an asymptotically de Sitter space  (see, e.g., \cite{Baskin,Vasy_2010}) and to involve  geometry into the analysis of the operators on the de~Sitter space (see, e.g., \cite{Birrell,Moschella,Parker-Toms,Tolman}).
\smallskip

In more explicit form the kernels can be written as follows 
\begin{eqnarray*}
K_0(z,t;M)
&  := &
4 ^ {-M}  e^{ t M}\big((1+e^{-t })^2 - z^2\big)^{  M    } \frac{1}{ [(1-e^{ -t} )^2 -  z^2]\sqrt{(1+e^{-t } )^2 - z^2} } \\
&   &
\times  \Bigg[  \big(  e^{-t} -1 +M(e^{ -2t} -      1 -  z^2) \big)
F \Big(\frac{1}{2}-M   ,\frac{1}{2}-M  ;1; \frac{ ( 1-e^{-t })^2 -z^2 }{( 1+e^{-t })^2 -z^2 }\Big)\nonumber  \\
&  &
\hspace{0.3cm}  +   \big( 1-e^{-2 t}+  z^2 \big)\Big( \frac{1}{2}+M\Big)
F \Big(-\frac{1}{2}-M   ,\frac{1}{2}-M  ;1; \frac{ ( 1-e^{-t })^2 -z^2 }{( 1+e^{-t })^2 -z^2 }\Big) \Bigg]\nonumber ,\\
K_1(z,t;M)
& :=  &
  4 ^{-M} e^{ Mt }  \big((1+e^{-t })^2 -   z  ^2\big)^{-\frac{1}{2}+M    } F\left(\frac{1}{2}-M   ,\frac{1}{2}-M  ;1;
\frac{ ( 1-e^{-t })^2 -z^2 }{( 1+e^{-t })^2 -z^2 } \right) \,.  
\end{eqnarray*}

We give  one example  of the equation  with the variable coefficients that is amenable to the integral transform method.
If $\Omega =\Pi $ is a non-Euclidean space of constant negative curvature  and the equation of the problems (\ref{1.22}) and (\ref{1.23})
is a non-Euclidean wave equation, then the explicit 
representation formulas are known (see, e.g., \cite{Helgason,L-P}) and the Huygens' principle is a consequence of those formulas.    
Thus, for a non-Euclidean wave equation, due to Theorem~\ref{T1.2}, the functions $v_f(x,t;b) $ and $v_\varphi (x,t) $ have  explicit 
representations, and the arguments of \cite{Yag_Galst_CMP,JMP2013}    allow us to derive for the solution $u(x,t)$ of the problem (\ref{1.20}) in the  de~Sitter   metric with   hyperbolic spatial geometry the explicit 
representation, the $L^p-L^q$ estimates, and to examine   the Huygens' principle. 
\medskip

 Thus, according to Theorem~\ref{T1.2} for the solution $\psi $ of the equation
\begin{equation}
\label{1.7a}
\label{1.9}
\psi _{tt} +   n  \psi _t - e^{-2 t} A(x,\partial_x)\psi  + m^2\psi =  f\,,
\end{equation}
due to the relation $u = e^{\frac{n}{2}t}\psi $, we obtain
\begin{eqnarray}
\label{1.10n}
\label{1.29l}
\psi   (x,t)
&  =  &
2  e^{-\frac{n}{2}t} \int_{ 0}^{t} db
  \int_{ 0}^{ e^{-b}- e^{-t}} dr  \,  e^{\frac{n}{2}b}v_f(x,r ;b) E(r,t; 0,b;M)  ,
\end{eqnarray}
where the function
$v_f(x,t;b)$
is defined by (\ref{1.22}), and
\[
M^2=\frac{n^2}{4}-m^2\,.
\]

Then, for the solution $\psi $ of the Cauchy problem
\begin{eqnarray}
\label{1.10old}
 \psi _{tt} +   n  \psi _t - e^{-2 t} A(x,\partial_x) \psi  + m^2\psi  =  0\,, \quad \psi  (x,0)= \psi _0 (x)\, , \quad \psi _t(x,0)=\psi _1 (x)\,,
\end{eqnarray}
due to the relation $u = e^{\frac{n}{2}t}\psi $, we obtain
\begin{eqnarray}
\label{1.11}
&  &
\psi  (x,t) \\
& = &
e^{-\frac{n-1}{2}t} v_{\psi _0}  (x, \phi (t))
+ \,  e^{-\frac{n}{2}t}\int_{ 0}^{1} v_{\psi _0}  (x, \phi (t)s)\Big(2  K_0(\phi (t)s,t;M)+ nK_1(\phi (t)s,t;M)\Big)\phi (t)\,  ds  \nonumber \\
& &
+\, 2e^{-\frac{n}{2}t}\int_{0}^1   v_{\psi  _1 } (x, \phi (t) s)
  K_1(\phi (t)s,t;M) \phi (t)\, ds
, \,\, t>0\,. \nonumber 
\end{eqnarray}

We stress here that the integral transform of this section allows to shrink (contract) the     time interval $[0,\infty)$ to the bounded interval $[0,1]$.

Then, for the linear strictly hyperbolic equations with smooth coefficients the integral representations hold also for the functions in the Sobolev and Besov spaces.  

Henceforth,  we use the following classification of the physical mass $m$. The mass $m$ is {\it large} if $m\geq n$, {\it small} if $m < n$, and {\it critical} (for the Huygens' principle) if   $m =\sqrt{n^2-1}/2 $.

\section{Large Mass. $B^{s,q}_{p}-B^{s',q}_{p'}$ Estimates  }
\label{S2}
\setcounter{equation}{0}

The following estimates proved by Brenner~\cite{Brenner1979} are crucial   for the estimates of this and the next sections. 
Let $\varphi _j= \varphi (2^{-j} \xi )  $, $j>0$, and $\varphi _0=1-\sum_{j=1}^\infty\varphi_j  $, where $\varphi \in C_0^\infty ({\mathbb R}^n) $ 
with $\varphi \geq 0 $  and \,\,supp$\,\varphi \subseteq \{ \xi \in {\mathbb R}^n\,;\, 1/2 <|\xi|< 2\}$, is that 
\[
\sum_{-\infty}^\infty \varphi  (2^{-j} \xi ) =1, \quad \xi \not=0\,.
\]
The norm $\|g\|_{B^{s,q}_{p}} $ of the Besov space $ B^{s,q}_{p}$ is defined as follows
\begin{eqnarray*}
&   &
\|v\|_{B^{s,q}_{p}}= \left(\sum_{j=0}^\infty \left( 2^{js} \| {\mathcal F}^{-1}\left( \varphi _j\hat{v}  \right) \|_{p} \right)^q \right) ^{1/q}\,,
\end{eqnarray*} 
where $\hat{v}  $ is the Fourier transform of $v$.

\begin{theorem} \mbox{\rm (Brenner~\cite{Brenner1979})}
\label{TBrenner}
Let $A=A(x,D)$ be a second order negative elliptic differential operator with real $C^\infty$-coefficients such that 
$ A(x,D) =A(\infty,D)$ for $|x|$ large enough. Let $u(t) = G_0(t)g_0+ G_1(t)g_1$ be the solution of 
\begin{eqnarray}
\label{2.1}
&   &
\partial_t^2 u - A(x,D)u =0, \quad x \in {\mathbb R}^n, \quad t \geq 0, \\
\label{2.2}
&  &
u (x,0) = g_0 (x),\quad   u_t (x,0)= g_1(x),\quad x \in {\mathbb R}^n\,.
\end{eqnarray}
Then for each $T<\infty$ there is a constant $C=C(T)$ such that if \, $(n+1)\delta \leq \nu +s-s' $,
\begin{eqnarray}
\label{Blplq}
&   &
\| G_\nu (t) g\| _{B^{s',q}_{p'}} \leq C(T)t^{\nu +s-s'-2n\delta } \|g\|_{B^{s,q}_{p}}, \quad 0 < t \leq T  \,.
\end{eqnarray}
Here $s,s' \geq 0$, $q\geq 1 $, $1\leq p\leq 2$, $1/p+1/p'=1$, and $\delta =1/p-1/2 $. 
\end{theorem}

The estimate (\ref{Blplq}) can be also written as
\begin{eqnarray*} 
&   &
\| G_\nu (t) g\| _{B^{s',q}_{p'}} \leq C(T)t^{\nu +s-s'-n\left(1/p-1/p' \right) } \|g\|_{B^{s,q}_{p}}, \quad 0 < t \leq T, \quad \nu =0,1 \,.
\end{eqnarray*}
We note that in this theorem $T<\infty$, but on the other hand, due to the integral transform of the previous section, it is possible to reduce    the problem with infinite time to the problem with the finite time, and apply Theorem~\ref{TBrenner} with $T=1$ only. 
We note that the condition $ A(x,D) =A(\infty,D)$ for $|x|$ large enough is not restrictive since the de~Sitter space-time has permanently bounded domain of influence. 

In fact, Theorem~3.2~\cite{Brenner1979} states similar result for the case of operator $ A=A(x,t,D) $ that makes possible to generalize some results of the present paper.

\subsection{Large Mass. $B^{s,q}_{p}-B^{s',q}_{p'}$  Estimates for  Equation without Source}
\label{S13}

The  decay estimates for the energy of the solution of the Cauchy problem for the 
wave equation without source can be proved 
by the representation formula, $L^1-L^\infty$ and  
$L^2-L^2$ estimates, and interpolation argument. (See, e.g., \cite[Theorem~2.1]{Racke}.) 
The proof of Theorem~\ref{TBrenner} given in   \cite{Brenner1979}  is based on the Fourier operators theory, 
microlocal consideration and dyadic decomposition of the phase space.   
In order to prove in the next theorem the similar result in the de~Sitter space-time,  we 
appeal to the representation formula provided by Theorem~\ref{T1.2} and then aplly Theorem~\ref{TBrenner}.

\begin{theorem} 
\label{T13.2}
Let $u(t) = G_{0,dS}(t)\varphi _0+ G_{1,dS}(t)\varphi _1$ be the solution of the Cauchy problem (\ref{1.20}) with $f=0$. Then  the operators 
$ G_{0,dS}(t)$ and $G_{1,dS}(t) $ satisfy the following  estimates
\begin{eqnarray*}
\|G_{0,dS}(t) \varphi_0 \|_{B^{s',q}_{p'}}  
&  \leq  & 
C_M  (1+ t  )^{1- \sgn M}(1- e^{-t}) ^{s-s'-2n\delta }  e ^{\frac{t}{2}}  \| \varphi_0  \|_{B^{s,q}_{p}}, 
\end{eqnarray*}
and 
\begin{eqnarray*}
\| G_{1,dS}(t)\varphi _1 \|_{B^{s',q}_{p'}}  
&  \leq  & 
C_M  (1+ t  )^{1- \sgn M}(1- e^{-t}) ^{ 1+  s-s'-2n\delta }   \|\varphi_1  \|_{B^{s,q}_{p}},  
\end{eqnarray*}
for all $t \in (0,\infty)$,  provided that  $(n+1)\delta \leq s-s' $, $1<p \leq 2$, $\frac{1}{p}+ \frac{1}{p'}=1$, $s-s'-2n\delta >-1 $,
and $\delta =1/p-1/2 $. 
 
They can be written as follows
\begin{eqnarray*}
\| G_{\nu, dS} (t) \varphi   \|_{B^{s',q}_{p'}}  
&  \leq  & 
C_M  e ^{\frac{t}{2}(1-\nu )}  (1+ t  )^{1- \sgn M}(1- e^{-t}) ^{\nu + s-s'-2n\delta }   \|\varphi   \|_{B^{s,q}_{p}}, \quad \nu =0,1 \,.
\end{eqnarray*}

\end{theorem}
\medskip

\noindent
{\it Proof.} 
We start with the operator $G_{1,dS}(t)$. Due  to Theorem~\ref{TBrenner}  for the solution $u=u (x,t)$ of the Cauchy 
problem (\ref{2.1})-(\ref{2.2}) with  $\varphi _0=0$ and according to 
(\ref{Blplq})
we have for $T=1$:

\begin{eqnarray*}
&   &
\| G_0 (t) g\| _{B^{s',q}_{p'}} \leq C(T)t^{ s-s'-2n\delta } \|g\|_{B^{s,q}_{p}}, \quad 0 < t \leq T ,  \quad (n+1)\delta \leq s-s' \,.
\end{eqnarray*}
Then, due to Theorem~\ref{T1.2} 
\begin{eqnarray*} 
\|  u(x,t) \|_{B^{s',q}_{p'}}  
& \leq &
2 \int_{0}^{1-e^{-t}}   
  \|v_{\varphi_1 } (x,r)  K_1(r,t;-iM)  \|_{B^{s',q}_{p'}}  \, dr \\ 
& \leq &
C \|\varphi_1  \|_{B^{s,q}_{p}}\int_{0}^{1-e^{-t}}   
 r^{s-s'-2n\delta }  \left| K_1(r,t;-iM)\right|   \, dr \\
& \leq &
C_M \|\varphi_1  \|_{B^{s,q}_{p}} e ^{ -t(s-s'-2n\delta) } \int_{0}^{e^{t } -1}   
y^{s-s'-2n\delta } \big( (e^{t } +1)^2 -y^2 \big)^{-\frac{1}{2}} \\
&  &
\hspace{2cm} \times 
\left( F\left(\frac{1}{2} ,\frac{1}{2}  ;1; 
\frac{ (e^{t } -1)^2 -y^2 }{(e^{t } +1)^2 -y^2 } \right) \right)^{1- \sgn M}\, dy \,.  
\end{eqnarray*}
To continue we prove the next simple generalization of    Lemma~9.2~\cite{Yag_Galst_CMP}. 
\begin{lemma}
\label{9.2CMP}
Assume that $0 \ge a>-1 $. Then 
\[
  \int_{ 0}^{ z- 1}  r^{ a }  \frac{1}{\sqrt{(z  + 1)^2    - r^2}}
F\left(\frac{1}{2},\frac{1}{2};1; 
\frac{  (z  - 1)^2  -r^2 }
{  (z  + 1)^2  -r^2  }  \right)  \, dr  
 \leq   
C     z^{-1}(z-1)^{ 1+  a } (1+ \ln z ),
\]
for all $z>1$.
\end{lemma}
\medskip

\noindent
{\it Proof.}   If $1 <  z \leq M$ with some constant $M$, then the argument $ $ of the hypergeometric functions is bounded, 
\begin{equation}
\label{argument}
\frac{ (z-1)^2 -r^2   }{ (z+1)^2 -r^2 } \leq \frac{ (z-1)^2  }{ (z+1)^2 } \leq \frac{ (M-1)^2  }{ (M+1)^2 }  <1 \quad \mbox{for all} \quad r \in (0,z-1),
\end{equation}
and
\begin{eqnarray*}
&  &
  \int_{ 0}^{ z- 1}  r^{ a }  \frac{1}{\sqrt{(z  + 1)^2    - r^2}}
F\left(\frac{1}{2},\frac{1}{2};1; 
\frac{  (z  - 1)^2  -r^2 }
{  (z  + 1)^2  -r^2  }  \right)  \, dr \\
 &  \leq &
C _M     (z-1)^{1+   a } , \quad \mbox{\rm for all } \quad 1 <  z \leq M\,.
\end{eqnarray*}
Hence, we can restrict ourselves to the case of large $z$, that is $z\geq M$. In particular, we choose $M>6$ and split integral into two parts:
\begin{eqnarray*}
&  &
  \int_{ 0}^{ z- 1}  r^{ a}  \frac{1}{\sqrt{(z  + 1)^2    - r^2}}
F\left(\frac{1}{2},\frac{1}{2};1; 
\frac{  (z  - 1)^2  -r^2 }
{  (z  + 1)^2  -r^2  }  \right)  \, dr \\
 & = &
  \int_{ 0}^{ \sqrt{(z  + 1)^2-8z}}  r^{ a }  \frac{1}{\sqrt{(z  + 1)^2    - r^2}}
F\left(\frac{1}{2},\frac{1}{2};1; 
1- \frac{  4z  }
{ (z  + 1)^2 -r ^2} \right)  \, dr \\
&     &
+   \int_{ \sqrt{(z  + 1)^2-8z}}^{ z- 1}  r^{ a}  \frac{1}{\sqrt{(z  + 1)^2    - r^2}}
F\left(\frac{1}{2},\frac{1}{2};1; 
1- \frac{  4z  }
{ (z  + 1)^2 -r ^2}\right)  \, dr \,.
\end{eqnarray*}
For the second part we have \,$z\ge M >6$\, and \,
$
r \ge \sqrt{(z  + 1)^2-8z} \,,
$
then
\begin{equation}
\label{9.2}
\frac{  4z  } { (z  + 1)^2 -r ^2}  \ge \frac{1}{2}    \Longrightarrow 0< 1- \frac{  4z  } { (z  + 1)^2 -r ^2} \le \frac{1}{2} 
\end{equation}
for such $r$ and $z$  implies
\begin{equation}
\label{9.3}
\left| F\left(\frac{1}{2},\frac{1}{2};1; 1-
\frac{  4z  }
{ (z  + 1)^2 -r ^2} \right) \right| \le C  \,.         
\end{equation}
Hence,   (\ref{9.2}) and (\ref{9.3}) imply 
\begin{eqnarray*}
&  &
\int_{\sqrt{(z  + 1)^2-8z}  }^{ z-1}  r^{ a }   \frac{1}{\sqrt{(z  + 1)^2  -r^2}}  
F\left(\frac{1}{2},\frac{1}{2};1; 1- \frac{  4z  }{ (z  + 1)^2 -r ^2} \right) \, dr  \\
&  \le &
C \int_{0}^{ z-1}  r^{ a }   \frac{1}{\sqrt{(z  + 1)^2  -r^2}}   \, dr \\
& \le  &
C  (1+z)^{a}   \quad \mbox{\rm for all} \quad z\ge M >6\,.
\end{eqnarray*}
For the first integral \, 
$ r \le \sqrt{(z  + 1)^2-8z} $  \, and  \, $ z\ge M >6 $ \,  imply  \,  $8z \leq (z  + 1)^2 - r^2 $. It follows  
\[
\left| F\left(\frac{1}{2},\frac{1}{2};1; 1-
\frac{  4z  }
{ (z  + 1)^2 -r ^2} \Bigg) \right| \le C \left| \ln \Bigg( \frac{  4z  }
{ (z  + 1)^2 -r ^2} \right)  \right|  \le C   (1+ \ln   z )  \,.       
\]
Then we obtain
\begin{eqnarray*}
&  &
  \int_{ 0}^{ \sqrt{(z  + 1)^2-8z}}  r^{ a }  \frac{1}{\sqrt{(z  + 1)^2    - r^2}}
F\left(\frac{1}{2},\frac{1}{2};1; 
1- \frac{  4z  }
{ (z  + 1)^2 -r ^2} \right)  \, dr \\
&   \leq   &
 C   (1+ \ln   z ) \int_{ 0}^{z-1}  r^{ a }  \frac{1}{\sqrt{(z  + 1)^2    - r^2}} \, dr \\
&   \leq   &
 C   (1+ \ln   z )  (1+z)^{a}   \,.
\end{eqnarray*}
Lemma is proved.
\hfill $\square$

In order to complete the proof of the theorem for this case, 
with $z= e^{t }$ and $s-s'-2n\delta >-1 $  we conclude  
\begin{eqnarray*}
&  &
\|  u(x,t) \|_{B^{s',q}_{p'}} \\
& \leq &
C_M \|\varphi_1  \|_{B^{s,q}_{p}} e ^{ -t(s-s'-2n\delta) } \\
&  &
\times \int_{0}^{e^{t } -1}   
y^{s-s'-2n\delta } \big( (e^{t } +1)^2 -y^2 \big)^{-\frac{1}{2}} 
\left( F\left(\frac{1}{2} ,\frac{1}{2}  ;1; 
\frac{ (e^{t } -1)^2 -y^2 }{(e^{t } +1)^2 -y^2 } \right) \right)^{1- \sgn M}\, dy \\
& \leq &
C_M \|\varphi_1  \|_{B^{s,q}_{p}} e ^{ -t(s-s'-2n\delta) }  e^{-t }(e^{t } -1)^{1+s-s'-2n\delta } (1+ t)^{1- \sgn M} \,.  
\end{eqnarray*}
Thus, in the case of $\varphi _0=0$ the theorem is proved. 
\medskip

Next we turn to the operator $G_{0,dS}(t)$.    Due  to Theorem~\ref{TBrenner} for the solution $u=u (x,t)$ of the Cauchy 
problem (\ref{2.1})-(\ref{2.2}) with  $\varphi _1=0$ and to Theorem~\ref{T1.2}  
we have: 
\begin{eqnarray*} 
\|  u(x,t) \|_{B^{s',q}_{p'}}
& \leq &
e^{t/2}\| v_{\varphi_0 } (x,1-e^{-t})\|_{B^{s',q}_{p'}}+2 \int_{0}^{1-e^{-t}}   
  \|v_{\varphi_0 } (x,r)  K_0(r,t;-iM)  \|_{B^{s',q}_{p'}}  \, dr \\  
& \leq &
e^{t/2}\| v_{\varphi_0 } (x,1-e^{-t})\|_{B^{s',q}_{p'}}+C \|\varphi_0  \|_{B^{s,q}_{p}}\int_{0}^{1-e^{-t}}   
 r^{s-s'-2n\delta }  \left| K_0(r,t;-iM)\right|   \, dr \\  
& \leq & 
C \Bigg( e ^{\frac{t}{2}} (1- e^{-t}) ^{s-s'-2n\delta } 
+ \int_{ 0}^{1- e^{-t}}  r^{s-s'-2n\delta }
|K_0(r,t;-iM)|  \,  dr \Bigg) \| \varphi_0  (x) \|_{B^{s,q}_{p}}.
\end{eqnarray*}
One can estimate  the last integral 
\begin{eqnarray*} 
&  &
\int_{ 0}^{1- e^{-t}}  r^{s-s'-2n\delta }
|K_0(r,t;-iM)|  \,  dr \\
& \leq & 
e^{ -t(s-s'-2n\delta ) }\int_{ 0}^{e^{t }-1}  y^{s-s'-2n\delta }
\frac{1}{ [(e^{t }-1)^2 -  y^2]\sqrt{(e^{t }+1)^2 - y^2} }\\
&   &
\times   \Bigg|  \big(  e^{t} - e^{2t} - iM(1 -     e^{ 2t} -  y^2) \big) 
F \Big(\frac{1}{2}+iM   ,\frac{1}{2}+iM  ;1; \frac{ ( e^{t }-1)^2 -y^2 }{( e^{t }+1)^2 -y^2 }\Big) \\
&  &
\hspace{1cm}  +   \big(  e^{2t} -1+  y^2 \big)\Big( \frac{1}{2}-iM\Big)
F \Big(-\frac{1}{2}+iM   ,\frac{1}{2}+iM  ;1; \frac{ ( e^{t }-1)^2 -y^2 }{( e^{t }+1)^2 -y^2 }\Big) \Bigg|  \,  d y \,.
\end{eqnarray*}
The following proposition with $a=s-s'-2n\delta $\, gives the remaining   estimate for that integral
and allow us to complete  the proof of the theorem. 
\begin{proposition}
\label{P13.4}
If  $a>-1 $, then
\begin{eqnarray*} 
&  &
\int_{ 0}^{z-1}  y^{a }
\frac{1}{ [(z-1)^2 -  y^2]\sqrt{(z+1)^2 - y^2} }\\
&   &
\times   \Bigg|  \big( z - z^{2 } - iM(1 -     z^{ 2 } -  y^2) \big) 
F \Big(\frac{1}{2}+iM   ,\frac{1}{2}+iM  ;1; \frac{ ( z-1)^2 -y^2 }{( z+1)^2 -y^2 }\Big) \\
&  &
\hspace{1cm}  +   \big(  z^{2 } -1+  y^2 \big)\Big( \frac{1}{2}-iM\Big)
F \Big(-\frac{1}{2}+iM   ,\frac{1}{2}+iM  ;1; \frac{ ( z-1)^2 -y^2 }{( z+1)^2 -y^2 }\Big) \Bigg|  \,  d y \\
& \leq &
 C_M z^{ -\frac{1}{2}}(z-1)^{1+ a }  \left(1+ \ln  z     \right)^{1- \sgn M}    \quad \mbox{ for all} \quad z>1. 
\end{eqnarray*}
\end{proposition}
\medskip

\noindent
{\it Proof.} (Comp. with  Prop.10.2~\cite{Yag_Galst_CMP}.) We follow the arguments have been used in the proof of Proposition~8.3~\cite{Yag_Galst_CMP}. If $1\leq  z \leq N$ with some constant $N$, then the argument  of the hypergeometric functions is bounded (\ref{argument}), 
and  the integral can be estimated by:
\begin{eqnarray*} 
&  &
\int_{ 0}^{z-1}  y^{a}
\frac{1}{ [(z-1)^2 -  y^2]\sqrt{(z+1)^2 - y^2} }\\
&   &
\times   \Bigg|  \big( z - z^{2 } - iM(1 -     z^{ 2 } -  y^2) \big) 
F \Big(\frac{1}{2}+iM   ,\frac{1}{2}+iM  ;1; \frac{ ( z-1)^2 -y^2 }{( z+1)^2 -y^2 }\Big) \\
&  &
\hspace{1cm}  +   \big(  z^{2 } -1+  y^2 \big)\Big( \frac{1}{2}-iM\Big)
F \Big(-\frac{1}{2}+iM   ,\frac{1}{2}+iM  ;1; \frac{ ( z-1)^2 -y^2 }{( z+1)^2 -y^2 }\Big) \Bigg|  \,  d y \\
& \leq &
C_M \int_{ 0}^{z-1}  y^{a}
\Bigg[ \frac{1}{  \sqrt{(z+1)^2 - y^2} }   
\Bigg\{ 1+ z^2\frac{1 }{(z +1)^2 -y^2 }   \Bigg\} 
 \Bigg]  \,  d y \\
& \leq &
C_Mz^{-1 }(z-1)^{ 1+a }  \quad \mbox{\rm for all}  \quad z \in [1,N]\,. 
\end{eqnarray*}
Thus, we can restrict ourselves to the  case of large $   z \geq N$ in both zones $Z_1(\varepsilon, z) $ and $ Z_2(\varepsilon, z )$, 
which are defined by
\begin{eqnarray*} 
Z_1(\varepsilon, z) 
& := &
\left\{ (z,r) \,\Big|\, \frac{ (z-1)^2 -r^2   }{ (z+1)^2 -r^2 } \leq \varepsilon,\,\, 0 \leq r \leq z-1 \right\} ,\\  
Z_2(\varepsilon, z) 
& := &
\left\{ (z,r) \,\Big|\, \varepsilon \leq  \frac{ (z-1)^2 -r^2   }{ (z+1)^2 -r^2 },\,\, 0 \leq r \leq z-1  \right\},
\end{eqnarray*}
respectively. In the first zone we have (7.11)~\cite{Yag_Galst_CMP}:
\begin{eqnarray*}
\hspace{-0.5cm} &   &
\Bigg|  \big( z -z^{ 2 }  - iM(1 -    z^{ 2 }-  r^2)    \big) 
F \Big(\frac{1}{2}+iM   ,\frac{1}{2}+iM  ;1; \frac{ ( z-1)^2 -r^2 }{( z+1)^2 -r^2 }\Big) \nonumber \\
\hspace{-0.5cm} &  &
\hspace{1cm}   +    \big(z^{2 }- 1+  r^2  \big) \Big( \frac{1}{2}-iM\Big)
F \Big(-\frac{1}{2}+iM   ,\frac{1}{2}+iM  ;1; 
\frac{ ( z-1)^2 -r^2 }{( z+1)^2 -r^2 }\Big)  \Bigg|  \nonumber \\
\hspace{-0.5cm} & \leq  &
\frac{1}{2}\big[(z-1)^2 -r^2 \big] +\frac{1}{8}\frac{ (z-1)^2 -r^2   }{ (z+1)^2 -r^2 }\nonumber \\
\hspace{-0.5cm} &  &
\times \left| (1-2 i M)  (-1+4 M^2 )  (r^2+z^2-1  )+2 (1+2 i M)^2  (-z^2+z +i M  (r^2+z^2 -1) )\right|  \nonumber  \\
\hspace{-0.5cm} &  & 
 + \Big( \left|  z -z^{ 2 }  - iM(1 -    z^{ 2 }-  r^2)      \right| 
+\left|  z^{2 }- 1+  r^2  \right|\Big)O\left(\left( \frac{ (z-1)^2 -r^2   }{ (z+1)^2 -r^2 }\right)^2\right) .
\end{eqnarray*} 
Consider therefore the following three estimates. For the first one we have (comp. with  Prop.10.2\cite{Yag_Galst_CMP})
\begin{eqnarray*}
A_{8}
& := &
\int_{ (z,r) \in Z_1(\varepsilon, z)  }  r^{a }  \frac{1}{     \sqrt{(z+1)^2-r^2}    }  \, dr  
  \leq   
C  z^{-1 }(z-1)^{ 1+a}  \quad \mbox{\rm for all}  \quad z \in [N,\infty)  \,.
\end{eqnarray*}
For $0 \geq a>-1 $ and $z\geq N$ the following integral can be easily estimated: 
\begin{eqnarray*} 
\int_{  0 }^{z-1}  r^{ a}    \frac{1}{ ((z+1)^2-r^2 )^{3/2}   }  dr 
& = &
\int_{  0 }^{z/2}  r^{ a}    \frac{1}{ ((z+1)^2-r^2 )^{3/2}   }  dr  + \int_{  z/2 }^{z-1}  r^{ a}    \frac{1}{ ((z+1)^2-r^2 )^{3/2}   }  dr  \\
& \leq  &
\frac{16}{9} z^{  -3}  \int_{  0 }^{z/2}  r^{ a }   dr  + \frac{z^{ a}}{4^{ a}} \int_{  z/2 }^{z-1}     \frac{1}{ ((z+1)^2-r^2 )^{3/2}   }  dr  \\
& \leq  &
C  z^{ a -3/2}  \quad \mbox{\rm for all}  \quad z \in [N,\infty)  \,.
\end{eqnarray*} 
Hence,
\begin{eqnarray*}
A_{9}
& := &
z^{2}\int_{  (z,r) \in Z_1(\varepsilon, z) } r^{ a }  
 \frac{1}{   \sqrt{(z+1)^2-r^2}    } \frac{ 1}{ (z+1)^2 -r^2 }  
 dr \\
& \leq   &
z^{2}\int_{ 0}^{z-1} r^{a}  
 \frac{1}{   \sqrt{(z+1)^2-r^2}    } \frac{ 1}{ (z+1)^2 -r^2 }  
 dr \\
& \leq   &
 C  z^{ -\frac{1}{2}}(z-1)^{1+ a}  
  \quad \mbox{\rm for all}  \quad z \in [N,\infty)\,,
\end{eqnarray*} 
and 
\begin{eqnarray*}
A_{10}
& := &
z^2\int_{  (z,r) \in Z_1(\varepsilon, z) }   r^{a} \frac{1}{ ((z-1)^2-r^2 ) \sqrt{(z+1)^2-r^2}    } \left(\frac{ (z-1)^2 -r^2   }{ (z+1)^2 -r^2 } \right)^2  dr \\
& \leq  &
z^2\int_{  (z,r) \in Z_1(\varepsilon, z) }  r^{ a}  \frac{1}{  \sqrt{(z+1)^2-r^2}    } \frac{ 1   }{  (z+1)^2 -r^2  }   dr \\
& \leq  &
 C  z^{ -\frac{1}{2}}(z-1)^{1+ a }  
  \quad \mbox{\rm for all}  \quad z \in [N,\infty) \,.
\end{eqnarray*}
Finally,
\begin{eqnarray*}
&  &
\int_{  (z,y) \in Z_1(\varepsilon, z) } \, y^{a}
\frac{1}{ [(z-1)^2 -  y^2]\sqrt{(z+1)^2 - y^2} }\\
&   &
\times   \Bigg|  \big( z - z^{2 } - iM(1 -     z^{ 2 } -  y^2) \big) 
F \Big(\frac{1}{2}+iM   ,\frac{1}{2}+iM  ;1; \frac{ ( z-1)^2 -y^2 }{( z+1)^2 -y^2 }\Big) \\
&  &
\hspace{1cm}  +   \big(  z^{2 } -1+  y^2 \big)\Big( \frac{1}{2}-iM\Big)
F \Big(-\frac{1}{2}+iM   ,\frac{1}{2}+iM  ;1; \frac{ ( z-1)^2 -y^2 }{( z+1)^2 -y^2 }\Big) \Bigg|  \,  d y \\& \leq  &
 C  z^{ -\frac{1}{2}}(z-1)^{1+ a}  
  \quad \mbox{\rm for all}  \quad z \in [1,\infty) \,.
\end{eqnarray*}
In the second zone we use (7.12)-(7.14)\cite{Yag_Galst_CMP} :
\begin{equation}
\label{10.4n}
\varepsilon \leq  \frac{ (z-1)^2 -r^2   }{ (z+1)^2 -r^2 } \leq 1 \quad \mbox{\rm and}  \quad 
\frac{ 1  }{ (z-1)^2 -r^2 }  \leq  \frac{ 1   }{ \varepsilon[(z+1)^2 -r^2] }\,.
\end{equation}
According to (7.3)~\cite{Yag_Galst_CMP}   
\begin{eqnarray*} 
\left| F\Big(\frac{1}{2}+iM   ,\frac{1}{2}+iM  ;1; 
z \Big) \right|
&  \leq  &
 C_{M }   \left(1 - \ln  (1-z)     \right) ^{1- \sgn M} \quad \mbox{\rm for all}\quad z \in [0 ,1)  \,.
\end{eqnarray*} 
Thus, the hypergeometric functions obey the estimates
\begin{equation}
\label{FGH+xat1}
\hspace{-0.5cm} \left| F\Big(-\frac{1}{2}+iM,\frac{1}{2}+iM;1;z   \Big) \right|  , \,\,   
\left| F\Big(\frac{1}{2}+iM,\frac{1}{2}+iM;1; z \Big) \right|  \leq C_M   \,\, \, \mbox{\rm for all}\,\,  z \in [\varepsilon ,1)  ,
\end{equation}
and
\begin{equation}
\label{10.8}
\left| F\Big(\frac{1}{2},\frac{1}{2};1; \frac{ (z-1)^2 -r^2   }{ (z+1)^2 -r^2 }    \Big) \right|   
  \leq   
C \left(1+ \ln  z     \right) 
, \quad \mbox{\rm for all}  \quad (z,r) \in Z_2(\varepsilon, z)  . 
\end{equation}
In the second zone we use (\ref{10.4n}), (\ref{FGH+xat1}),   and  (\ref{10.8}).  Thus, we have to estimate the next two integrals:
\begin{eqnarray*} 
A_{11} 
& :=   &
z^{2}\int_{  (z,r) \in Z_2(\varepsilon, z) }   r^{ a } \frac{1}{     ((z-1)^2-r^2 ) \sqrt{(z+1)^2-r^2}    }  \,dr \,, \\
A_{12} 
& :=   & 
z^{2}\left(1+ \ln  z     \right)^{1- \sgn M} \int_{  (z,r) \in Z_2(\varepsilon, z) }   r^{ a }  
\frac{1}{     ((z-1)^2-r^2 ) \sqrt{(z+1)^2-r^2}    } \,dr \,.
\end{eqnarray*}
We apply (\ref{10.4n}) to $A_{11}$ and obtain
\begin{eqnarray*} 
A_{11} 
 \leq   
C_\varepsilon z^{2}\int_{  (z,r) \in Z_2(\varepsilon, z) }   r^{a}   
\frac{ 1   }{ [(z+1)^2 -r^2] }\frac{1}{    \sqrt{(z+1)^2-r^2}    }  \, dr
 \leq  
C_\varepsilon z^{ -\frac{1}{2}}(z-1)^{1+ a}    
\end{eqnarray*} 
for all $ z \in [1,\infty)$, while
\begin{eqnarray*} 
A_{12} 
& \leq    & 
C_\varepsilon z^{ -\frac{1}{2}}(z-1)^{1+ a}  \left(1+ \ln  z     \right)^{1- \sgn M}   
\quad \mbox{\rm for all}  \quad z \in [1,\infty)  \,.
\end{eqnarray*} 
The proposition is proved.
\hfill $\square$
\smallskip

To complete the proof of the theorem we write 
\begin{eqnarray*} 
&  &
\int_{ 0}^{1- e^{-t}}  r^{s-s'-2n\delta }
|K_0(r,t;-iM)|  \,  dr \\
& \leq & 
e^{ -t[s-s'-2n\delta ] }\int_{ 0}^{e^{t }-1}  y^{s-s'-2n\delta  }
\frac{1}{ [(e^{t }-1)^2 -  y^2]\sqrt{(e^{t }+1)^2 - y^2} }\\
&   &
\times   \Bigg|  \big(  e^{t} - e^{2t} - iM(1 -     e^{ 2t} -  y^2) \big) 
F \Big(\frac{1}{2}+iM   ,\frac{1}{2}+iM  ;1; \frac{ ( e^{t }-1)^2 -y^2 }{( e^{t }+1)^2 -y^2 }\Big) \\
&  &
\hspace{1cm}  +   \big(  e^{2t} -1+  y^2 \big)\Big( \frac{1}{2}-iM\Big)
F \Big(-\frac{1}{2}+iM   ,\frac{1}{2}+iM  ;1; \frac{ ( e^{t }-1)^2 -y^2 }{( e^{t }+1)^2 -y^2 }\Big) \Bigg|  \,  d y \\
& \leq & 
Ce^{ -t[\frac{1}{2} + s-s'-2n\delta ] } (e^t-1)^{1+ s-s'-2n\delta }  (1+ t  )^{1- \sgn M}\,.
\end{eqnarray*}
Thus,
\begin{eqnarray*} 
&  &
\|  u(x,t) \|_{B^{s',q}_{p'}} \\
& \leq &
C \Bigg( e ^{\frac{t}{2}} (1- e^{-t}) ^{s-s'-2n\delta } 
+ \int_{ 0}^{1- e^{-t}}  r^{s-s'-2n\delta }
|K_0(r,t;-iM)|  \,  dr \Bigg) \| \varphi_0  (x) \|_{B^{s,q}_{p}}\\
& \leq &
C \Bigg( e ^{\frac{t}{2}} (1- e^{-t}) ^{s-s'-2n\delta } 
+ e^{ -t[\frac{1}{2} + s-s'-2n\delta ] } (e^t-1)^{1+ s-s'-2n\delta  }  (1+ t  )^{1- \sgn M}\Big)  \Bigg) \| \varphi_0  (x) \|_{B^{s,q}_{p}}\\
& \leq &
C  (1+ t  )^{1- \sgn M}e ^{\frac{t}{2}} (1- e^{-t}) ^{s-s'-2n\delta }  \| \varphi_0  (x) \|_{B^{s,q}_{p}}.
\end{eqnarray*}
Theorem is proved. \hfill $\square$
\medskip

\subsection{Large Mass. $B^{s,q}_{p}-B^{s',q}_{p'}$ Estimates for the Equation with  Source}

In general, the  Duhamel's principle allows us to reduce the case with a  source term to the case of the Cauchy problem without
source term and consequently to derive the decay estimates for the equation. For the equation under consideration 
we   appeal to the representation formula of Theorem~\ref{T1.2}. 
In this subsection we consider the Cauchy problem  for the equation with the source term and with zero initial data.   

\begin{theorem} 
\label{T2.5}
Let $u=u(x,t)$ be solution of the Cauchy problem  
\begin{equation}
\label{1.20b}
\cases{
u_{tt} - e^{-2t}A(x,\partial_x) u - M^2 u= f, \quad  x \in {\mathbb R}^n \,,\,\, t >0,\cr
  u(x,0)= 0\, , \quad u_t(x,0)=0,\quad x \in {\mathbb R}^n \,.
}
\end{equation}
Then for $n \geq 2$ one has the following  estimate
\begin{eqnarray*}
\hspace{-0.5cm} &  &
\|  u(x,t) \|_{B^{s',q}_{p'}}\\
\hspace{-0.5cm} &  & \le 
C_M \int_{ 0}^{t} db\, \|  f(x, b)  \|_{B^{s,q}_{p}}  
  \int_{ 0}^{ e^{-b}- e^{-t}} r^{s-s'-2n\delta  }  \frac{\left(F\left(\frac{1}{2},\frac{1}{2};1; 
\frac{ (e^{-b}- e^{-t})^2-r^2}
{  (e^{-b}+ e^{-t})^2-r^2}   \right) \right)^{1- \sgn M} }{\sqrt{(e^{-t}  + e^{-b})^2    - r^2}}
dr\, 
\end{eqnarray*}
for all $t>0$, provided that  $1<p \leq 2$, $\frac{1}{p}+ \frac{1}{p'}=1$, $   s-s'-2n\delta>-1$,  $s,s'\geq  0$, $(n+1)\delta \leq s-s'$, and $\delta =1/p-1/2$.
\end{theorem}
\medskip

\noindent
{\it Proof.} We use the representation  given by Theorem~\ref{T1.2}. 
Due to Theorem~\ref{TBrenner} for the   equation of (\ref{1.20b}), we have
\begin{eqnarray*} 
\|  u(x,t) \|_{B^{s',q}_{p'}}
&  \leq   &
2   \int_{ 0}^{t} db
  \int_{ 0}^{\phi (t)- \phi (b)} \| v_f(x,r ;b)\|_{B^{s',q}_{p'}} |E(r,t;0,b;-iM)|  \, dr \,,    
\end{eqnarray*}
where $\phi (t):= 1-e^{-t} $, and, consequently, 
\begin{eqnarray*}
\hspace{-0.5cm} &  &
\|  u(x,t) \|_{B^{s',q}_{p'}}\\
\hspace{-0.5cm} & &  \le 
C_M\int_{ 0}^{t} db
  \int_{ 0}^{ e^{-b}- e^{-t}}   \| v_f(x,r ;b)\|_{B^{s',q}_{p'}}   \frac{\left(F\left(\frac{1}{2},\frac{1}{2};1; 
\frac{ (e^{-b}- e^{-t})^2-r^2}
{  (e^{-b}+ e^{-t})^2-r^2}  \right)  \right)^{1- \sgn M}}{\sqrt{(e^{-t}  + e^{-b})^2    - r^2}}
dr \\ 
\hspace{-0.5cm} &  & \le 
C_M \int_{ 0}^{t} db\, \|  f(x, b)  \|_{B^{s,q}_{p}}  
  \int_{ 0}^{ e^{-b}- e^{-t}} r^{  s-s'-2n\delta  }  \frac{\left(F\left(\frac{1}{2},\frac{1}{2};1; 
\frac{ (e^{-b}- e^{-t})^2-r^2}
{  (e^{-b}+ e^{-t})^2-r^2}   \right) \right)^{1- \sgn M} }{\sqrt{(e^{-t}  + e^{-b})^2    - r^2}}
dr\,.
\end{eqnarray*}
Theorem is proved. \hfill $\square$
\smallskip

We are going to transform the estimate of the last theorem to a more compact form. 
To this aim we  estimate  for  $s-s'-2n\delta>-1 $ the  last integral of the right hand side. 
If we replace $e^{-(b -t)} $ with $z := e^{-(b -t)}>  1$, then the integral will be simplified:
\begin{eqnarray*}
&  &
\int_{ 0}^{ e^{-b}- e^{-t}} r^{ s-s'-2n\delta }  \frac{1}{\sqrt{(e^{-t}  + e^{-b})^2    - r^2}}
\left( F\left(\frac{1}{2},\frac{1}{2};1; 
\frac{ (e^{-b}- e^{-t})^2-r^2}
{  (e^{-b}+ e^{-t})^2-r^2}   \right)  \right)^{1- \sgn M} dr \\
& = &
e^{-t(s-s'-2n\delta)}\int_{ 0}^{ z- 1} y^{ s-s'-2n\delta }  \frac{1}{\sqrt{(z  + 1)^2    - y^2}}
\left(F\left(\frac{1}{2},\frac{1}{2};1; 
\frac{ ( z  - 1 )^2-y^2}
{   (z  + 1) ^2-y^2}   \right)  \right)^{1- \sgn M} dy
\end{eqnarray*}
Then we apply Lemma~\ref{9.2CMP} and obtain the following corollary. 

\begin{corollary}
\label{C12.3}
Let $u=u(x,t)$ be a solution of the Cauchy problem (\ref{1.20b}). Then for $n\geq 2$ one has the following estimate  
\[
\|  u(x,t) \|_{B^{s',q}_{p'}}  \le 
C_M \int_{ 0}^{t} db\, \|  f(x, b)  \|_{B^{s,q}_{p}}  
  e^{b}\left( e^{-b}- e^{-t}  \right)^{ 1+  s-s'-2n\delta   } \left(1+  t-  b    \right)^{1- \sgn M}\,  db\,,
\] 
provided that  $1<p \leq 2$, $\frac{1}{p}+ \frac{1}{p'}=1$, $s,s'\geq 0$,  $s-s'-2n\delta>-1 $, $(n+1)\delta \leq s-s' $.
\end{corollary}
\medskip
 
\noindent
{\it Proof.}  
Indeed, we apply Lemma~\ref{9.2CMP}  with $z= e^{t -b}$ to the right-hand side of the estimate given by Theorem~\ref{T2.5} :
\begin{eqnarray*}
\hspace{-0.5cm} &  &
\|  u(x,t) \|_{B^{s',q}_{p'}}\\
\hspace{-0.5cm} &  & \le 
C_M \int_{ 0}^{t} db\, \|  f(x, b)  \|_{B^{s,q}_{p}}  
  \int_{ 0}^{ e^{-b}- e^{-t}} r^{ s-s'-2n\delta  }  \frac{\left(F\left(\frac{1}{2},\frac{1}{2};1; 
\frac{ (e^{-b}- e^{-t})^2-r^2}
{  (e^{-b}+ e^{-t})^2-r^2}   \right) \right)^{1- \sgn M} }{\sqrt{(e^{-t}  + e^{-b})^2    - r^2}}
dr\\
\hspace{-0.5cm} &  & \le 
C_M \int_{ 0}^{t} db\, \|  f(x, b)  \|_{B^{s,q}_{p}}  
 e^{-t(  s-s'-2n\delta )}    z^{-1}(z-1)^{ 1+   s-s'-2n\delta  } (1+ \ln z )^{1- \sgn M}\\
\hspace{-0.5cm} &  & \le 
C_M \int_{ 0}^{t}   \|  f(x, b)  \|_{B^{s,q}_{p}}  
  e^{b}\left( e^{-b}- e^{-t}  \right)^{ 1+  s-s'-2n\delta   } \left(1+  t-  b    \right)^{1- \sgn M}\,  db\,.
\end{eqnarray*} 
The corollary is proved.
\hfill $\square$

\subsection{Large mass. $B^{s,q}_{p}-B^{s',q}_{p'}$ estimates for the covariant equation}
\label{SS2.3}

Then, for the solution $\psi $ of the Cauchy problem
\begin{equation}
\label{1.10}
\psi _{tt} +   n   \psi _t - e^{-2 t} A(x,\partial_x) \psi  + m^2\psi  =  0\,, \quad \psi  (x,0)= \psi _0 (x)\, , \quad \psi  _t(x,0)=\psi _1 (x)\,,
\end{equation}
due to the relation $u = e^{\frac{n}{2}t}\psi $, we obtain (\ref{1.11}). 

\noindent
{\bf Estimates for  equation with  source.}
Let $u=u(x,t)$ be a solution of the Cauchy problem (\ref{1.20b}). Then according to Corollary~\ref{C12.3},
for the solution $\psi  $ of the equation (\ref{1.7a}),  
due to (\ref{1.10n}), we obtain
\begin{equation}
\label{2.4new}
\|\psi  (x,t) \|_{B^{s',q}_{p'}}  
\le   
C_Me^{-\frac{n}{2}t}\int_{ 0}^{t}e^{\frac{n}{2}b} \|  f(x, b)  \|_{B^{s,q}_{p}}  
 e^{b}\left( e^{-b}- e^{-t}  \right)^{ 1+   s-s'-2n\delta   } \left(1+  t-  b    \right)^{1- \sgn M}\,  db .
\end{equation}
For $M>0$   we obtain 
\[
\| \psi  (x,t) \|_{B^{s',q}_{p'}}  
\le   
C_Me^{-\frac{n}{2}t}\int_{ 0}^{t}e^{\frac{n}{2}b} \|  f(x, b)  \|_{B^{s,q}_{p}}  
 e^{b}\left( e^{-b}- e^{-t}  \right)^{ 1+  s-s'-2n\delta  } \,  db \,,
\]
while for $M=0$   we obtain 
\[
\| \psi (x,t) \|_{B^{s',q}_{p'}}  
\le   
C e^{-\frac{n}{2}t}\int_{ 0}^{t}e^{\frac{n}{2}b} \|  f(x, b)  \|_{B^{s,q}_{p}}  
 e^{b}\left( e^{-b}- e^{-t}  \right)^{ 1+ s-s'-2n\delta  } \left(1+  t-  b    \right) \,  db .
\]
In particular, for $s=s' $, $\delta =0 $ and $p=p'=2$, that is for   the Sobolev space $H_{(s)}({\mathbb R}^n)$    we have 
\begin{equation}
\label{2.13}
\|\psi  (x,t) \|_{ H_{(s)}({\mathbb R}^n)}  
\le   
C_Me^{-\frac{n}{2}t}\int_{ 0}^{t}e^{\frac{n}{2}b} \|  f(x, b)  \|_{H_{(s)}({\mathbb R}^n)}   \left(1+  t-  b    \right)^{1- \sgn M}
 \,  db. 
\end{equation}
Here the rates of exponential factors are independent of the curved mass ${\mathcal M} $ if ${\mathcal M}\not=0$ and, consequently, of the mass $m$.
\medskip

\noindent 
 {\bf  Equation without source.} 
According to Theorem~\ref{T13.2} the solution  $u=u(x,t)$ of the Cauchy problem (\ref{1.20}) with $f=0$  satisfies the following   estimate   
\begin{eqnarray*}
\|u(x,t) \|_{B^{s',q}_{p'}}  
&  \leq  & 
C_M  (1+ t  )^{1- \sgn M}(1- e^{-t}) ^{ s-s'-2n\delta }  e ^{\frac{t}{2}}  \| \varphi_0  (x) \|_{B^{s,q}_{p}}, 
\end{eqnarray*}
if $ (n+1)\delta \leq  s-s'$, and 
\begin{eqnarray*}
\|  u(x,t) \|_{B^{s',q}_{p'}}  
&  \leq  & 
C_M  (1+ t  )^{1- \sgn M}(1- e^{-t}) ^{1+    s-s'-2n\delta }   \|\varphi_1  \|_{B^{s,q}_{p}},  
\end{eqnarray*}
if $(n+1)\delta \leq s-s'$ and $s-s'-2n\delta >-1 $  

Thus, for the solution $\psi  $ of the Cauchy problem (\ref{1.10}),  
due to the relation $u = e^{\frac{n}{2}t}\psi  $, we obtain the decay  estimate 
\begin{equation}
\label{1.13}
\label{2.4}
\|  \psi  (x,t) \|_{B^{s',q}_{p'}} 
  \leq   
C_M e^{-\frac{n}{2}t} (1+ t  )^{1- \sgn M}(1- e^{-t}) ^{  s-s'-2n\delta}\Big\{ e ^{\frac{t}{2}}  \| \psi _0  (x) \|_{B^{s,q}_{p}}
+ (1- e^{-t}) \|\psi _1  \|_{B^{s,q}_{p}} 
 \Big\} 
\end{equation}
for all $t>0$,  while 
\[ 
\| \psi  (x,t) \|_{B^{s',q}_{p'}} 
  \leq   
C_M e^{-\frac{n}{2}t}(1+ t  )^{1- \sgn M} \Big\{ e ^{\frac{t}{2}}  \| \psi _0  (x) \|_{B^{s,q}_{p}}
+  \|\psi _1  \|_{B^{s,q}_{p}} 
 \Big\}  
\]
for large $t$. Here the rate of decay is independent of the curved mass ${\mathcal M}$ if ${\mathcal M}\not=0$ and, consequently, of the mass $m$.

\section{Small Mass.  $B^{s,q}_{p}-B^{s',q}_{p'}$ Estimates}
\label{S3}

\setcounter{equation}{0}

{\bf Equation without source.}    We derive the estimates of solutions of the covariant Klein-Gordon equation from the estimates of solutions 
of the non-covariant equation with the imaginary mass.  
\begin{theorem} 
\label{T13.2s}
\label{T3.1}
Let  $\psi  =\psi  (x,t)$ be a  solution of the Cauchy problem 
\[
 \psi  _{tt} +   n  \psi  _t - e^{-2 t} A(x,\partial_x) \psi   \pm m^2\psi   =  0\,, \quad \psi  (x,0)= \psi _0 (x)\, , \quad \psi  _t(x,0)=\psi _1 (x)\,.
\] 
Define $M=\sqrt{\frac{n^2}{4}-m^2}>0$    for the case of ``plus'', and  $M=\sqrt{\frac{n^2}{4}+m^2}$ for the case of ``minus''.
Then, in the case of ``plus'' with  $m \in (0, \sqrt{n^2-1}/{2})$ and in the case of ``minus'' for all $m>0$,  it 
satisfies the following $B^{s,q}_{p}-B^{s',q}_{p'}$ estimate
\[
\| \psi (x,t) \|_{B^{s',q}_{p'}} 
  \leq   
C_{m,n,p,s}   (1- e^{-t}) ^{  s-s'-2n\delta }e^{( M -\frac{n}{2})t}\Big\{   \| \psi _0   \|_{B^{s,q}_{p}}
+ (1- e^{-t}) \|\psi _1  \|_{B^{s,q}_{p}} 
 \Big\} 
\]
for all $t \in (0,\infty)$, provided that  $s,s' \geq 0$, $(n+1)\delta \leq s-s'$, $\delta =1/p-1/2 $, and $1<p \leq 2$, $\frac{1}{p}+ \frac{1}{p'}=1$. 
The above estimate   holds also if $m \in ( \sqrt{n^2-1}/{2},n/2)$ and $\psi _0=0 $ that is,
\[
\| \psi (x,t) \|_{B^{s',q}_{p'}} 
  \leq   
C_{m,n,p,s}   (1- e^{-t}) ^{ 1+ s-s'-2n\delta }e^{( M -\frac{n}{2})t}      (1- e^{-t}) \|\psi _1  \|_{B^{s,q}_{p}}  \,.
\] 
\end{theorem}
\medskip

\noindent
{\it Proof.} First we consider the case of $\psi _1=0 $. Then
\begin{eqnarray*}  
\psi  (x,t)  
& = &
e^{-\frac{n-1}{2}t} v_{\psi _0}  (x, \phi (t))\\
&  &
+ \,  e^{-\frac{n}{2}t}\int_{ 0}^{1} v_{\psi _0}  (x, \phi (t)\tau )\big(2  K_0(\phi (t)\tau ,t;M)+ nK_1(\phi (t)\tau ,t;M)\big)\phi (t)\,  d\tau 
\end{eqnarray*}
and, consequently, 
\begin{eqnarray} 
\label{2.9} 
\|\psi  (x,t) \|_{B^{s',q}_{p'}}  
&   \leq &  
e^{-\frac{n-1}{2}t} \|  v_{\psi _0}  (x, \phi (t)) \|_{B^{s',q}_{p'}}  \\
&  &
+ \,  e^{-\frac{n}{2}t}\int_{ 0}^{1} \|  v_{\psi _0}  (x, \phi (t)\tau )\|_{B^{s',q}_{p'}}  \big|2  K_0(\phi (t)\tau ,t;M)+ nK_1(\phi (t)\tau ,t;M)\big|\phi (t)\,  d\tau \,. \nonumber
\end{eqnarray}
If $n \geq 2$, then for the solution $v = v (x,t)$ of the Cauchy problem (\ref{2.1})-(\ref{2.2}) 
with $\psi  (x) \in C_0^\infty({\mathbb R}^n)$  we apply Theorem~\ref{TBrenner}
\begin{eqnarray*} 
&  & 
\| v_\psi   (x ,t) \|_{B^{s',q}_{p'}}
\le C 
t^{s-s'-2n\delta }\|\psi   \|_{B^{s,q}_{p}}   \quad \mbox{\rm for all} \quad  t >0,
\end{eqnarray*}
provided that \, $s, s'\ge 0$, \, $(n+1)\delta \leq s-s'$,  $1<p\le 2$, $\frac{1}{p}+\frac{1}{p'}=1$. 
Hence,
\[ 
 \| v_{\psi _0}  (x, \phi (t))\|_{B^{s',q}_{p'}}
 \leq  
C\phi (t)^{s-s'-2n\delta }\|\psi  _0 \|_{B^{s,q}_{p}}  \quad \mbox{\rm for all} \quad  t >0,
\]
where $\phi (t)= 1-e^{-t}$. Consequently, for the first term of the right-hand side of (\ref{2.9}) we have
\begin{eqnarray*}  
e^{-\frac{n-1}{2}t} \|  v_{\psi _0}  (x, \phi (t)) \|_{B^{s',q}_{p'}}  
& \leq  &
C  e^{-\frac{n-1}{2}t}(1-e^{-t})^{s-s'-2n\delta }\|\psi  _0 \|_{B^{s,q}_{p}}  \quad \mbox{\rm for all} \,\,  t >0,
\end{eqnarray*}
while for the second  term we obtain
\begin{eqnarray*} 
&  &
e^{-\frac{n}{2}t}\int_{ 0}^{1} \|  v_{\psi _0}  (x, \phi (t)\tau )\|_{B^{s',q}_{p'}} \big|2  K_0(\phi (t)\tau ,t;M)+ nK_1(\phi (t)\tau ,t;M)\big|\phi (t)\,  d\tau \\
& \leq &
 \|\psi _0 \|_{B^{s,q}_{p}} e^{-\frac{n}{2}t} \int_{ 0}^{1} \phi (t)^{s-s'-2n\delta} \tau ^{s-s'-2n\delta}
 \left( \big|2  K_0(\phi (t)\tau ,t;M)\big|+ n\big|K_1(\phi (t)\tau ,t;M)\big|\right) \phi (t)\,  d\tau  \,.
\end{eqnarray*}
We have to estimate the following two integrals of the last inequality:
\begin{eqnarray*} 
&  & 
 \int_{ 0}^{1} \phi (t)^{s-s'-2n\delta} \tau ^{s-s'-2n\delta}
\big|  K_0(\phi (t)\tau ,t;M)\big|   \phi (t)\,  d\tau  
\end{eqnarray*}
and 
\begin{eqnarray*} 
&  & 
 \int_{ 0}^{1} \phi (t)^{s-s'-2n\delta} \tau ^{s-s'-2n\delta}
 \big|K_1(\phi (t)\tau ,t;M)\big| \phi (t)\,  d\tau  \,,
\end{eqnarray*}
where  $\phi (t)= 1-e^{-t}$ and $t>0$. We   apply  Lemma~2.3~\cite{JMAA2012} and Lemma~2.4~\cite{JMAA2012} (it is important here that $M>1/2 $) in the case of $a=s-s'-2n\delta$ and arrive at the following estimates 
\begin{eqnarray*} 
&  &
  \int_{ 0}^{1} \phi (t)^{s-s'-2n\delta} \tau ^{s-s'-2n\delta}
\big|  K_0(\phi (t)\tau ,t;M)\big|   \phi (t)\,  d\tau  \nonumber \\ 
& \leq &
C_{M,n,p,q,s}     e^ {-Mt-(s-s'-2n\delta)t}(e^t-1)^{1 +s-s'-2n\delta} (e^t+1)^{2 M-1}  \quad \mbox{\rm for all}  \quad t>0 \,,  
\end{eqnarray*}
and for $M>1/2$, while 
\begin{eqnarray*} 
&  & 
\int_{ 0}^{1} \phi (t)^{s-s'-2n\delta} \tau ^{s-s'-2n\delta}
\big|  K_1(\phi (t)\tau ,t;M)\big|   \phi (t)\,  d\tau    \\
& \leq &
C_{M,n,p,q,s}  e^{-Mt -t(s-s'-2n\delta)}(e^{t }-1)^{1+s-s'-2n\delta} (e^{t }+1)^{2 M-1}  \quad \mbox{\rm for all}  \quad t>0\,,
\end{eqnarray*}
and for $M>0$. 

Now consider the case of $\psi _0=0 $. We have 
\begin{eqnarray*} 
\psi  (x,t) 
& = & \, 2e^{-\frac{n}{2}t}\int_{0}^1   v_{\psi  _1 } (x, \phi (t)\tau ) 
  K_1(\phi (t)\tau ,t;M) \phi (t)\, d\tau 
, \quad x \in {\mathbb R}^n, \,\, t>0\,.
\end{eqnarray*}
Then, according to Lemma~2.3~\cite{JMAA2012} we have
\begin{eqnarray*} 
\|\psi  (x,t) \|_{B^{s',q}_{p'}} 
& \leq  & \, 2e^{-\frac{n}{2}t}\int_{0}^1  \| v_{\psi  _1 } (x, \phi (t)\tau ) \|_{B^{s',q}_{p'}}
|  K_1(\phi (t)\tau ,t;M)| \phi (t)\, d\tau  \\
& \leq  & 
Ce^{-\frac{n}{2}t}\int_{0}^1  \| \psi  _1   \|_{B^{s,q}_{p}}\phi (t)^{s-s'-2n\delta} \tau ^{s-s'-2n\delta}
|  K_1(\phi (t)\tau ,t;M)| \phi (t)\, d\tau \\
& \leq  & 
C_Me^{-\frac{n}{2}t} \| \psi  _1   \|_{B^{s,q}_{p}}e^{-Mt-(s-s'-2n\delta)t} (e^t-1)^{1+s-s'-2n\delta} (e^t+1)^{2M-1} 
\end{eqnarray*}
that completes the proof of that case. Theorem is proved. \hfill $\square$
\bigskip

\noindent
{\bf $\bf B^{s,q}_{p}-B^{s',q}_{p'}$ Estimates for  Equation with  Source} 

\begin{theorem} 
\label{T11.1}
Let $\psi  =\psi  (x,t)$ be a solution of the Cauchy problem  
\begin{equation}
\label{2.5}
  \psi _{tt} +   n   \psi  _t - e^{-2 t} A(x,\partial_x) \psi   \pm m^2\psi  =  f\,, \quad \psi   (x,0)= 0\, , \quad \psi  _t(x,0)=0\,.
\end{equation}
Define $M=\sqrt{\frac{n^2}{4}-m^2} $  and    $m < n/2$   for the case of ``plus'', and  $M=\sqrt{\frac{n^2}{4}+m^2} $ for the case of ``minus''.
Then  $\psi  =\psi  (x,t)$ satisfies the following $B^{s,q}_{p}-B^{s',q}_{p'}$ estimate:
\begin{eqnarray} 
\label{3.3}
\|  \psi  (x ,t) \|_{ B^{s',q}_{p'}} 
&  \leq   &
 C_M e^{ -M t  } e^{-\frac{n}{2}t}  e^{-t(s-s'-2n\delta )}   \\
  &  &
\times 
\int_{ 0}^{t} e^{\frac{n}{2}b} e^{ M b  }(e^{t-b}-1)^{1+s-s'-2n\delta} (e^{t-b}+1)^{2 M-1}\| f(x ,b)  \|_{B^{s,q}_{p}}  \,db\,,  \nonumber 
\end{eqnarray} 
for all $t>0$, provided that \, $s, s'\ge 0$, $s,s' \geq 0$, $(n+1)\delta \leq s-s'$,  $1<p\le 2$, $\frac{1}{p}+\frac{1}{p'}=1$.
\end{theorem}
\medskip

\noindent
{\it Proof.} From (\ref{1.29l}) we have  
\begin{eqnarray*} 
\psi (x,t) 
&  =  &
2  e^{-\frac{n}{2}t} \int_{ 0}^{t} db
  \int_{ 0}^{ e^{-b}- e^{-t}} dr  \,  e^{\frac{n}{2}b} v_f(x,r ;b)  4 ^{-M}  e^{ M(b+t) } \Big((e^{-t }+e^{-b})^2 - r^2\Big)^{-\frac{1}{2}+M    } \\
  &  &
\times F\Big(\frac{1}{2}-M   ,\frac{1}{2}-M  ;1; 
\frac{ ( e^{-b}-e^{-t })^2 -r^2 }{( e^{-b}+e^{-t })^2 -r^2 } \Big) , 
\end{eqnarray*}
where according to (\ref{Blplq}) of Theorem~\ref{TBrenner} we can write 
\begin{eqnarray*} 
&  & 
\|v_f(x,r ;b)   \|_{ B^{s',q}_{p'}}
\le C 
r^{s-s'-2n\delta}\| f(x ,b)  \|_{B^{s,q}_{p}}  \quad \mbox{\rm for all} \quad  r >0\,.
\end{eqnarray*}
Hence,
\begin{eqnarray*} 
&  &
\| \psi (x ,t) \|_{ B^{s',q}_{p'}} \\
&  \leq   &
2  e^{-\frac{n}{2}t} \int_{ 0}^{t} db
  \int_{ 0}^{ e^{-b}- e^{-t}} dr  \,  e^{\frac{n}{2}b}\| v_f(x ,r ;b) \|_{ B^{s',q}_{p'}}  4 ^{-M}  e^{ M(b+t) } \Big((e^{-t }+e^{-b})^2 - r^2\Big)^{-\frac{1}{2}+M    } \\
  &  &
\times \left| F\Big(\frac{1}{2}-M   ,\frac{1}{2}-M  ;1; 
\frac{ ( e^{-b}-e^{-t })^2 -r^2 }{( e^{-b}+e^{-t })^2 -r^2 } \Big) \right|\\
&  \leq   &
 C_M e^{ M t  } e^{-\frac{n}{2}t} \int_{ 0}^{t} e^{\frac{n}{2}b} e^{ M b  }\| f(x ,b)  \|_{B^{s,q}_{p}}\,db \int_{ 0}^{ e^{-b}- e^{-t}} \, r^{s-s'-2n\delta}     \Big((e^{-t }+e^{-b})^2 - r^2\Big)^{-\frac{1}{2}+M    }
  \\
  &  &
\times 
  \left| F\Big(\frac{1}{2}-M   ,\frac{1}{2}-M  ;1; 
\frac{ ( e^{-b}-e^{-t })^2 -r^2 }{( e^{-b}+e^{-t })^2 -r^2 } \Big) \right| dr \,.  
\end{eqnarray*} 
We set   $  r=ye^{-t}$ and  obtain 
\begin{eqnarray*} 
&  &
\| \psi  (x ,t) \|_{ B^{s',q}_{p'}} \\
&  \leq   &
 C_M e^{ -M t  } e^{-\frac{n}{2}t}  e^{-t(s-s'-2n\delta) } \int_{ 0}^{t} e^{\frac{n}{2}b} e^{ M b  }\| f(x ,b)  \|_{B^{s,q}_{p}} \,db   \\
  &  &
 \times 
\int_{ 0}^{ e^{t-b}- 1} \, y^{s-s'-2n\delta}   \Big(( e^{t-b}+1)^2 - y^2\Big)^{-\frac{1}{2}+M    } \left| F\Big(\frac{1}{2}-M   ,\frac{1}{2}-M  ;1; 
\frac{ ( e^{t-b}-1)^2 -y^2 }{( e^{t-b}+1)^2 -y^2 } \Big) \right| dy \,.  
\end{eqnarray*}
Hence, we have to estimate the integral 
\begin{eqnarray*} 
\int_{ 0}^{ z- 1} \, y^{a}   \Big(( z+1)^2 - y^2\Big)^{-\frac{1}{2}+M    } \left| F\Big(\frac{1}{2}-M   ,\frac{1}{2}-M  ;1; 
\frac{ (z-1)^2 -y^2 }{(z+1)^2 -y^2 } \Big) \right| dy ,  
\end{eqnarray*}
where $z=e^{t-b} >1$ and $a= s-s'-2n\delta>-1$. On the other hand, since $M>0$,   we have   
\begin{eqnarray*} 
&  &
\int_{ 0}^{ z- 1} \, y^{a}   \Big(( z+1)^2 - y^2\Big)^{-\frac{1}{2}+M    } \left| F\Big(\frac{1}{2}-M   ,\frac{1}{2}-M  ;1; 
\frac{ (z-1)^2 -y^2 }{(z+1)^2 -y^2 } \Big) \right| dy \\
&  \leq   &
C_M\int_{ 0}^{ z- 1} \, y^{a}   \Big(( z+1)^2 - y^2\Big)^{-\frac{1}{2}+M    }  dy\\
&  =   &
C_M\frac{1}{1+a}(e^{t-b}-1)^{1+a} (e^{t-b}+1)^{-1+2 M} F\left(\frac{1+a}{2},\frac{1}{2}-M;\frac{3+a}{2};\frac{(e^{t-b}-1)^2}{(e^{t-b}+1)^2}\right)\\
& \leq    &
C_M\frac{1}{1+a}(e^{t-b}-1)^{1+a} (e^{t-b}+1)^{-1+2 M} \, ,
\end{eqnarray*} 
since
\begin{eqnarray*} 
\left| F\left(\frac{1+a}{2},\frac{1}{2}-M;\frac{3+a}{2};\frac{(e^{t-b}-1)^2}{(e^{t-b}+1)^2}\right) \right|
& \leq    & C_M \,  .
\end{eqnarray*} 
This completes the proof of the estimate (\ref{3.3}).
Theorem is proved. \hfill $\square$ 

The following corollaries can be easily verified. 
\begin{corollary} 
\label{C3.3}
\label{C2.7}
Let $\psi  =\psi  (x,t)$ be a solution of the Cauchy problem considered in Theorem~\ref{T11.1}. Then for $n \geq 2$  and $M>0 $ one has the following  estimate
\begin{eqnarray*} 
\| \psi (x ,t) \|_{ B^{s',q}_{p'}} 
&  \leq   &
 C_M   e^{-(\frac{n}{2}-M)t}    
\int_{ 0}^{t} e^{(\frac{n}{2}-M)b} e^{-b(s-s'-2n\delta)}  \| f(x ,b)  \|_{B^{s,q}_{p}} \,db ,  
\end{eqnarray*}
provided that  \, $s, s'\ge 0$,  $(n+1)\delta \leq s-s'$,  $1<p\le 2$, $\frac{1}{p}+\frac{1}{p'}=1$, $1<p \leq 2$, $s-s'-2n\delta>-1$, and $\delta =1/p-1/2 $.
\end{corollary}

\begin{corollary} 
\label{C3.4}
Let $\psi  =\psi  (x,t)$ be a solution of the Cauchy problem considered in Theorem~\ref{T11.1}. Then for $n \geq 2$  and $M>0 $ one has the following  estimate
\begin{eqnarray*} 
\| \psi (x ,t) \|_{ H_{(s)} ({\mathbb R}^n)} 
&  \leq   &
 C_M   e^{-(\frac{n}{2}-M)t}    
\int_{ 0}^{t} e^{(\frac{n}{2}-M)b}  \| f(x ,b)  \|_{ H_{(s)} ({\mathbb R}^n)} \,db ,  
\end{eqnarray*}
\end{corollary}

\section{Critical Mass $m^2=(n^2-1)/4$}
\label{S4}

\setcounter{equation}{0}

The value $m^2= (n^2-1)/4 $ of the mass $m$, when  $M=1/2$, dramatically simplifies the function $E(x,t;x_0,t_0;M) $, and
consequently, $K_0 (z,t;M)$ and $K_1 (z,t;M)$.   
Indeed, in that case we have
\begin{eqnarray*}
E\left(x,t;x_0,t_0;\frac{1}{2}\right) =
 \frac{1}{2} e^{ \frac{1}{2}(t_0+t) } ,\quad E\left(z,t;0,b;\frac{1}{2}\right) = \frac{1}{2} e^{ \frac{1}{2}(b+t) }\,,
\end{eqnarray*}
while
\begin{eqnarray*}
K_0\left(z,t;\frac{1}{2} \right)
  =
- \frac{1}{4}  e^{ \frac{1}{2}t }  ,\qquad
K_1\left(z,t;\frac{1}{2} \right)
   =   \frac{1}{2}  e^{ \frac{1}{2}t }   \,.
 \end{eqnarray*}

\subsection{Integral Transform and Huygens' Principle}
\label{SS3.1cr}
In \cite{JMP2013} for the equation with $A(x,\partial_x)=\Delta  $ is proved that for {\it critical value} $  (n^2-1)/4 $ of mass is the only value which makes equation Huygensian, that is, solutions of the equation obey Huygens' principle. For the equations with $x$-dependent coefficient    we have  the following  theorem.
\begin{theorem}
\label{T3}
Assume that $m= \sqrt{n^2-1}/2$. Then    the solutions of the equation
\begin{eqnarray}
\label{25}
  \psi  _{tt} +   n   \psi  _t - e^{-2 t} A(x,\partial_x) \psi   + m^2\psi  = 0,
\end{eqnarray}
 obey the strong Huygens' Principle,
whenever  the solutions of the   equation 
\[
 u_{tt}    -   A(x,\partial_x) u   = 0,
\]
obey it. 
\end{theorem}
\medskip

\noindent
{\it Proof.} 
For the solution (\ref{1.29l}) of the equation (\ref{1.7a}) with the source term it follows
\begin{eqnarray*}
\psi   (x,t)
&  =  &
    e^{-\frac{n-1}{2}t}\int_{ 0}^{t}  e^{\frac{n+1}{2}b} db
  \int_{ 0}^{ e^{-b}- e^{-t}}v(x,r ;b)  \,   dr    \,,
\end{eqnarray*}
where the function $v(x,r ;b) $ is defined by 
(\ref{1.22}). In fact,
if we denote by $V_{f}(x,t;b)$ the solution of the problem
\begin{eqnarray*}
V_{tt}-  A(x,\partial_x)  V =0, \quad V(x,0)= 0, \quad V_t(x,0)=f (x,b)\,,
\end{eqnarray*}
then
\[
v(x,t;b) =\frac{\partial }{\partial t}V_{f}(x,t;b)\,.
\]
Hence,
\begin{eqnarray}
\label{4.4}
\psi  (x,t)
&  =  &
  e^{-\frac{n-1}{2}t}\int_{ 0}^{t}  e^{\frac{n+1}{2}b}
   V_{f}(x,e^{-b}- e^{-t};b)   \,   db .
\end{eqnarray}
Further, due to (\ref{1.11}), for the solution $\psi  $ of (\ref{1.10old})    we have
\begin{eqnarray*}
\psi   (x,t)
& = &
e^{-\frac{n-1}{2}t} v_{\psi _0}  (x, 1-e^{-t})
+ \,   \frac{n-1}{2}e^{-\frac{n-1}{2}t}\int_{ 0}^{1-e^{-t}} v_{\psi _0}  (x, s )        \,  ds  \nonumber \\
& &
+\,  e^{-\frac{n-1}{2}t}\int_{0}^{1-e^{-t}}   v_{\psi _1 } (x, s )
     \, ds
, \quad x \in {\mathbb R}^n, \,\, t>0\,,
\end{eqnarray*}
where the functions $v_{\psi _0 }    $ and $v_{\psi _1 }  $ are
defined by (\ref{1.6}). Now, if we denote by $V_{\varphi }$ the
solution of the problem
\begin{eqnarray*}
V_{tt}-   A(x,\partial_x)  V =0, \quad V(x,0)= 0, \quad V_t(x,0)=\varphi (x),
\end{eqnarray*}
then
\[
v_{\varphi }(x,t) =\frac{\partial }{\partial t}V_{\varphi }(x,t)\,,
\]
and
\begin{eqnarray*}
\psi  (x,t)
& = &
e^{-\frac{n-1}{2}t} v_{\psi _0}  (x, 1-e^{-t})
+ \,   \frac{n-1}{2}e^{-\frac{n-1}{2}t}    V_{\psi _0}  (x, 1-e^{-t} )    \nonumber \\
& &
+\,  e^{-\frac{n-1}{2}t}  V_{\psi _1}  (x, 1-e^{-t} )
, \quad x \in {\mathbb R}^n, \,\, t>0\,    ,
\end{eqnarray*}
or, equivalently,
\begin{eqnarray*}
\psi (x,t)
& = &
e^{-\frac{n-1}{2}t}  \left(  \frac{\partial V_{\psi _0 }}{\partial t} \right) (x, 1-e^{-t})
+ \,   \frac{n-1}{2}e^{-\frac{n-1}{2}t}    V_{\psi _0}  (x, 1-e^{-t} )    \nonumber \\
& &
+\,  e^{-\frac{n-1}{2}t}  V_{\psi _1}  (x, 1-e^{-t} )
, \quad x \in {\mathbb R}^n, \,\, t>0\,.
\end{eqnarray*}
Theorem is proved. \hfill $\square$

\smallskip

\subsection{The Critical Case. Asymptotic Expansions of Solutions   at Infinite Time}
\label{SS3.2cr}

 In this subsection we generalize the results of \cite{JMP2013}, about the  asymptotic expansions of solutions  at the large time, to the equations with $x$-dependent coefficients.
For $\varphi  \in C_0^\infty ({\mathbb R}^n)$ let $v_{\varphi }   (x, t) $ be a 
 solution   of the Cauchy problem
\begin{eqnarray}
\label{4.6}
u_{tt}- A(x,\partial_x) u =0, \quad u(x,0)= 0, \quad u_t(x,0)=\varphi  (x)\,.
\end{eqnarray}
Denote
\begin{eqnarray*}
v_{\varphi }   (x):= v_{\varphi }   (x, 1)\,, \qquad V_{\varphi }(x):=V_{\varphi }(x, 1)\,.
\end{eqnarray*}
In order to write the complete asymptotic expansion of the
solutions, we define the functions
\begin{eqnarray*}
V_{\varphi }^{(k)}(x)
&: = &
\frac{(-1)^k}{k!}\left[\left( \frac{\partial }{\partial t} \right)^k V_{\varphi }(x, t) \right]_{t=1}
\in C_0^\infty ({\mathbb R}^n)\,,\quad k=1,2,\ldots\,.
\end{eqnarray*}
Then, for every integer  $N \geq 1$ we have
\begin{eqnarray*}
V_{\varphi }(x,1-e^{-t}) = \sum_{k=0}^{N-1} V_{\varphi }^{(k)}(x)e^{-kt} + R_{V_\varphi ,N}(x,t),  \quad
R_{V_\varphi ,N} \in   C ^\infty \,,
\end{eqnarray*}
where with the constant $ C(\varphi ) $ the remainder  $ R_{V_\varphi ,N} $ satisfies the inequality
\begin{eqnarray*}
| R_{V_\varphi ,N}(x,t) | \leq C(\varphi ) e^{-Nt}\quad  \,\, \mbox{\rm for all}\quad  x \in {\mathbb R}^n \quad  \mbox{\rm and all}\quad  t \in [0,\infty)  \,.
\end{eqnarray*}
Moreover, the support of the remainder $ R_{V_\varphi ,N} $ is in the cylinder
\begin{eqnarray*}
 \mbox{\rm supp}\,   R_{V_\varphi ,N}   \subseteq
   \{ x \in {\mathbb R}^n \,;\, \mbox{\rm dist} (x,   \mbox{\rm supp}\, \varphi ) \leq c_0\, \}\times [0,\infty)  \,.
\end{eqnarray*}
where $c_0$ is speed of propagation of the equation (\ref{4.6}). Analogously, we define 
\begin{eqnarray*}
v_{\varphi }^{(k)}(x)
& = &
\frac{(-1)^k}{k!}\left[\left( \frac{\partial }{\partial t} \right)^k v_{\varphi }(x, t) \right]_{t=1}
\in C_0^\infty ({\mathbb R}^n)\,,\quad k=1,2,\ldots\,,
\end{eqnarray*}
and the remainder $ R_{v_\varphi ,N} $
\begin{eqnarray*}
v_{\varphi }(x,1-e^{-t}) = \sum_{k=0}^{N-1} v_{\varphi }^{(k)}(x)e^{-kt} + R_{v_\varphi ,N}(x,t),  \quad
R_{v_\varphi ,N} \in   C ^\infty \,,
\end{eqnarray*}
such that
\begin{eqnarray*}
| R_{v_\varphi ,N}(x,t) | \leq C(\varphi ) e^{-Nt}\quad  \,\,
\mbox{\rm for all}\quad  x \in {\mathbb R}^n \quad  \mbox{\rm and
all}\quad  t \in [0,\infty).
\end{eqnarray*}
Further,  we introduce   the polynomial in $z $ with the smooth in $x \in {\mathbb R}^n$ coefficients  as follows:
\[
\psi _{asympt}^{(N)}  (x,z)
 =
 z^{\frac{n-1}{2}}  \left(    \sum_{k=0}^{N-1} v_{\varphi_0 }^{(k)}(x)z^k
+     \frac{n-1}{2}     \sum_{k=0}^{N-1} V_{\varphi_0 }^{(k)}(x)z^k  \right)
+     z^{\frac{n-1}{2}} \sum_{k=0}^{N-1} V_{\varphi_1 }^{(k)}(x)z^k      ,
\]
where $ x \in {\mathbb R}^n$, $z \in {\mathbb C} $. This allows us to    prove the following asymptotic expansion
\begin{eqnarray*}
\psi  (x,t) = \psi _{asympt}^{(N)}  (x,e^{-t})  + O(e^{-Nt-\frac{n-1}{2}t})
\end{eqnarray*}
for large $t$ uniformly for  $x\in {\mathbb R}^n$. Thus, we have proved the next theorem.

\begin{theorem}
\label{Tasymp}
Suppose that  $m= \sqrt{ n^2-1 }/2 $. Then, for every integer positive $N$   the solution of the equation (\ref{25}) with the initial values
$ \varphi_0, \varphi _1  \in C_0^\infty ({\mathbb R}^n) $   has the following asymptotic expansion at infinity:
\begin{eqnarray*}
\psi  (x,t) \sim  \psi _{asympt}^{(N)}  (x,e^{-t})\,,
\end{eqnarray*}
in the sense that  for every integer positive $N$  the following
estimate is valid:
\begin{eqnarray*}
\| \psi  (x,t) - \psi _{asympt}^{(N)}  (x,e^{-t}) \|_{L^\infty  ({\mathbb R}^n)}
& \leq  &
C(\varphi_0, \varphi _1) e^{-Nt-\frac{n-1}{2}t} \quad \,\, \mbox{  for large }\,\, t\,.
\end{eqnarray*}
\end{theorem}
\smallskip

\begin{remark}
If we take into account the relation
$v_{\varphi }   (x, t) = \frac{\partial }{\partial t} V_{\varphi } (x, t) $, then
\begin{eqnarray*}
&  &
v_{\varphi }^{(k)}   (x) = - (k+1)  V_{\varphi }^{(k+1)}   (x)
\end{eqnarray*}
and, consequently, the function $\Phi_{asympt}^{(N)}  (x,z) $ can be rewritten as follows:
\begin{eqnarray*}
&  &
\psi _{asympt}^{(N)}  (x,z) \\
& = &
 z^{\frac{n-1}{2}}  \left(    \sum_{k=0}^{N-1} (-1)(k+1) V_{\varphi_0 }^{(k+1)}(x)z^k
+     \frac{n-1}{2}     \sum_{k=0}^{N-1} V_{\varphi_0 }^{(k)}(x)z^k  \right)
+     z^{\frac{n-1}{2}} \sum_{k=0}^{N-1} V_{\varphi_1 }^{(k)}(x)z^k   \\
& = &
 z^{\frac{n-1}{2}}    \sum_{k=0}^{N-1}\left( \frac{n-1}{2}      V_{\varphi_0 }^{(k)}(x)  -(k+1) V_{\varphi_0 }^{(k+1)}(x)
+   V_{\varphi_1 }^{(k)}(x)\right) z^k .
\end{eqnarray*}
\end{remark}

\subsection{The Critical Case. $\bf B^{s,q}_{p}-B^{s',q}_{p'}$-Estimates}
\label{SS3.3cr}

\begin{lemma}
\label{L1} 
Suppose that  $m= \sqrt{ n^2-1  }/2$,  \, $(n+1)\delta \leq 1+ s-s'$, $s, s' \geq 0 $, $1\leq p \leq 2$, $1/p+1/p'=1$,  $s-s'-2n\delta >-2 $, and $\delta =1/p-1/2 $.  If $ \psi  _0
=\psi _1 =0$,  then for the solution $\psi 
= \psi  (x,t)$ of the equation (\ref{1.9}), (\ref{0.5})  the following estimate
holds
\begin{eqnarray*} 
\|\psi  (x,t) \|_{ B^{s',q}_{p'}}   
& \leq  &
C e^{-\frac{n-1}{2}t}\int_{ 0}^{t}  e^{\frac{n+1}{2}b}
(e^{-b}- e^{-t})^{1+ s-s'-2n\delta } \|
  f(x, b)   \|_{ B^{s,q}_{p}}  \,   db
, \quad  t>0\,.
\end{eqnarray*}
For the solution $\psi  = \psi  (x,t)$  of the  Cauchy problem
(\ref{1.10old}) if  $ \psi _0 =0$, and
  $s-s'-2n\delta >-1 $, $(n+1)\delta \leq 1+ s-s'$, then
\begin{eqnarray*}
\|  \psi   (x,t)  \|_{ B^{s',q}_{p'}}  
& \leq  &
C e^{-\frac{n-1}{2}t} (1-e^{-t})^{1+ s-s'-2n\delta }  \| \psi _1   \|_{ B^{s,q}_{p}}  
, \quad  t>0\,.
\end{eqnarray*}
If  $ \psi _1 =0$ and $(n+1)\delta \leq s-s'$,   then
\begin{eqnarray*}
\|\psi   (x,t)  \|_{ B^{s',q}_{p'}} 
& \leq  &
C e^{-\frac{n-1}{2}t} (1-e^{-t})^{  s-s'-2n\delta } \| \psi  _0   \|_{ B^{s,q}_{p}}
, \quad  t>0\,.
\end{eqnarray*} 
\end{lemma}
\smallskip

\noindent
{\it Proof.} According to  (\ref{4.4}) and Theorem~\ref{TBrenner}, if $\psi  _0 =\psi _1 =0$ and $(n+1)\delta \leq 1+ s-s'$, then
\begin{eqnarray*} 
\|\psi   (x,t)  \|_{ B^{s',q}_{p'}}  
& \leq  &
\| e^{-\frac{n-1}{2}t}\int_{ 0}^{t}  e^{\frac{n+1}{2}b}
   V_{f}(x,e^{-b}- e^{-t};b)   \,   db  \|_{ B^{s',q}_{p'}} \\
& \leq  &
 e^{-\frac{n-1}{2}t}\int_{ 0}^{t}  e^{\frac{n+1}{2}b}
 \| V_{f}(x,e^{-b}- e^{-t};b)   \|_{ B^{s',q}_{p'}} \,   db \\
& \leq  &
C e^{-\frac{n-1}{2}t}\int_{ 0}^{t}  e^{\frac{n+1}{2}b}
(e^{-b}- e^{-t})^{1+ s-s'-2n\delta }  \|
  f(x, b)    \|_{ B^{s,q}_{p}}\,   db
, \quad  t>0\,.
\end{eqnarray*}
In particular,
\begin{eqnarray*}
& &
\|\psi   (x,t)  \|_{ B^{s',q}_{p'}}\\
& \leq  &
C \left( \sup_{0\leq b\leq t}  \|
  f(x, b)   \|_{ B^{s,q}_{p}}\right)e^{-\frac{n-1}{2}t}\int_{ 0}^{t}  e^{\frac{n+1}{2}b}
(e^{-b}- e^{-t})^{1+ s-s'-2n\delta }  \,   db \\
& \leq  &
C \left( \sup_{0\leq b\leq t}  \|
  f(x, b)  \|_{ B^{s,q}_{p}}\right)e^{-\frac{n-1}{2}t}e^{-t( 1+ s-s'-2n\delta )} \int_{ 0}^{t}  e^{\frac{n+1}{2}b}
(e^{t-b}- 1)^{1+ s-s'-2n\delta }  \,   db
, \quad  t>0\,.
\end{eqnarray*}
For the case $s=s' $, $\delta =0 $, and $p=p'=2$ we obtain the estimate
\begin{eqnarray*}
\|\psi   (x,t) \|_{H_{(s)}({\mathbb R}^n)  }
& \leq  &
C \left( \sup_{0\leq b\leq t}  \|
  f(x, b)   \|_{H_{(s)}({\mathbb R}^n)  }\right)e^{-\frac{n-1}{2}t}\int_{ 0}^{t}  e^{\frac{n-1}{2}b}
  \,   db\\
& \leq  &
C \left( \sup_{0\leq b\leq t}  \|
  f(x, b)  \|_{H_{(s)}({\mathbb R}^n)  }\right),  \quad  t>0\,.
\end{eqnarray*}
 Further, if $ f \equiv 0$, $ \psi  _0 =0$, and $(n+1)\delta \leq 1+ s-s'$, then
\begin{eqnarray*}
\|\psi  (x,t) \|_{ B^{s',q}_{p'}}
& \leq  &
C e^{-\frac{n-1}{2}t} (1-e^{-t})^{1+s-s'-2n\delta}  \| \psi  _1  \|_{ B^{s,q}_{p}}
, \quad  t>0\,,
\end{eqnarray*}
If $ f \equiv 0$, $ \psi  _1 =0$, and  $(n+1)\delta \leq  s-s'$, then
\begin{eqnarray*}
\| \psi   (x,t) \|_{ B^{s',q}_{p'}}
& \leq  &
C e^{-\frac{n-1}{2}t} (1-e^{-t})^{s-s'-2n\delta} \| \psi  _0  \|_{ B^{s,q}_{p}}
, \quad  t>0\,.
\end{eqnarray*}
The lemma is proved. \hfill $ \square $ 
\smallskip

\section{Global Existence. Small Data Solutions}
\label{S5}

\setcounter{equation}{0}

Let $\psi  =\psi  (x,t)$ be a solution of the Cauchy problem (\ref{2.5})
with either $M=\sqrt{\frac{n^2}{4}-m^2} $ and  $m \in (0,n/2)$ for the case of ``plus'', or  $M=\sqrt{\frac{n^2}{4}+m^2} $ for the case of ``minus''. 
Then for $n \geq 2$,  according to Corollary~\ref{C3.3} one has the following  estimate 
\begin{eqnarray*} 
\| \psi (x ,t) \|_{ B^{s,q}_{2}} 
&  \leq   &
 C_M   e^{-(\frac{n}{2}-M)t}    
\int_{ 0}^{t} e^{(\frac{n}{2}-M)b}  \| f(x ,b)  \|_{B^{s,q}_{2}} \,db ,  
\end{eqnarray*}
For the equation with ``plus'' and large  mass,  $m \geq  n/2 $, and with the curved mass\\
$M=\sqrt{m^2 - n^2/4}$, according to (\ref{2.4new}) one has the following  estimate 
\[
\| \psi  (x,t) \|_{ B^{s,q}_{2}} 
\le   
C_Me^{-\frac{n}{2}t}\int_{ 0}^{t}e^{\frac{n}{2}b} (1+ t -b )^{1- \sgn M}\|  f(x, b)  \|_{B^{s,q}_{2}}   
 \,  db. 
\]
Here the rate of exponential factors  is independent of the curved mass $M$ and, consequently, of the mass $m$.
 These statements   follow  immediately from  Corollary~\ref{C2.7} and Section~\ref{SS2.3}.  
\medskip

The last estimates and the fixed point theorem allow us to prove global existence in the Cauchy problem for the semilinear equation 
\[
\psi _{tt} +   n  \psi  _t - e^{-2 t} A(x,\partial_x)\psi   +  m^2 \psi  =   F(\psi  )\,.
\]

We study the Cauchy problem   (\ref{NWE}), (\ref{ICPHI})
    through the integral equation.
To determine that integral equation we  appeal to the operator 
\[
G:={\mathcal K}\circ {\mathcal EE_A} \,,
\]
where 
\[
{\mathcal  EE_A} [f](x,t;b)= v(x,t;b)
\]
and the function 
$v(x,t;b)$   
is a solution to the Cauchy problem for the   equation (\ref{1.6}),  
while ${\mathcal K}$ is introduced either by (\ref{1.29l}),
\begin{eqnarray}
\label{operK}
{\mathcal K}[v]  (x,t) 
 & :=  &
2  e^{-\frac{n}{2}t} \int_{ 0}^{t} db
  \int_{ 0}^{ e^{-b}- e^{-t}} dr  \,  e^{\frac{n}{2}b}v(x,r ;b) E(r,t; 0,b;-iM) \\
 &  =  &
2  e^{-\frac{n}{2}t} \int_{ 0}^{t} db
  \int_{ 0}^{ e^{-b}- e^{-t}} dr  \,  e^{\frac{n}{2}b}{\mathcal  EE_A} [f](x,r ;b) E(r,t; 0,b;-iM) , \nonumber
\end{eqnarray}
for the large mass $m$, or by 
\begin{eqnarray} 
\label{operKsm}
{\mathcal K}[v]  (x,t) 
&  :=  &
2   e^{-\frac{n}{2}t}\int_{ 0}^{t} db
  \int_{ 0}^{ e^{-b}- e^{-t}} dr  \,  e^{\frac{n}{2}b} v(x,r ;b) E(r,t; 0,b;M) \\    
& = &
2   e^{-\frac{n}{2}t}\int_{ 0}^{t} db
  \int_{ 0}^{ e^{-b}- e^{-t}} dr  \,  e^{\frac{n}{2}b}\,{\mathcal EE_A} [f](x,r ;b) E(r,t; 0,b;M) \nonumber \,,
\end{eqnarray}
for the small mass $m$. 
Thus, the Cauchy problem (\ref{NWE}), (\ref{ICPHI}) leads to the following integral equation
\begin{eqnarray} 
\label{5.1}
\label{4.5}
\label{5.5}
\psi (x,t)
 = 
\psi _0(x,t) + 
G[ F(\psi  ) ] (x,t)    \,. 
\end{eqnarray}
Every solution to  the  equation (\ref{NWE}) solves also the last integral equation with some function $\psi  _0 (x,t)$, which, in fact, is a solution of the Cauchy problem   (\ref{1.10}).

 \subsection{Solvability of the Integral Equation associated with Klein-Gordon Equation}
\label{S7} 
\label{SS3.3}
\label{SS5.1}

\medskip

\noindent
Consider the integral equation
(\ref{5.1})
where $\psi  _0=\psi  _0 (x,t)$ is a given function. Every solution to   the equation (\ref{NWE}) 
solves also the last integral equation with some function $\psi  _0=\psi  _0 (x,t)$.
In order to solve the integral equation (\ref{5.1}), we  apply the Banach fixed-point theorem.  
To estimate nonlinear term we use the   Lipschitz condition (${\mathcal L}$).
Evidently, the condition ($\mathcal L$) imposes some restrictions on $n$, $\alpha $, $s$. 
Now we consider the integral equation  (\ref{5.1}),  
where the function $\psi _0\in C([0,\infty);B^{s,q}_{p})$ is given. 
 We note here that any classical  
solution to the  equation (\ref{NWE}) solves also the integral equation
(\ref{5.1}) with some function $\psi _0(t,x)$, which is  classical  
solution to the Cauchy problem for the linear equation  (\ref{1.10old}).
\medskip

Solvability of the integral equation (\ref{5.1}) is determined by the operator $G={\mathcal K}\circ {\mathcal EE_A} $. 
We start with the case of Sobolev space $ H_{(s)} ({\mathbb R}^n)$ with $ s > n/2$, which is an algebra.   
In the next theorem operator ${\mathcal K} $ is generated by the linear part of the equation (\ref{NWE}). 
\begin{theorem}
\label{TIE}  
Assume that  $F$  is Lipschitz continuous with exponent $\alpha >0 $ in the  space $H_{(s)} ({\mathbb R}^n)$, $ s > n/2$, and $F(x,0)=0$ for all $x \in {\mathbb R}^n$.
Then, there exists sufficiently small $\varepsilon _0>0 $ such that, for every given function $ \psi _0(x ,t) \in X({\varepsilon ,s,\gamma_0}) $, $\varepsilon <\varepsilon _0 $, such that  
\begin{eqnarray*}  
&  &
\sup_{t \in [0,\infty)}  e^{\gamma_0 t}  \|\psi  _0(\cdot ,t) \|_{H_{(s)} ({\mathbb R}^n)}  < \varepsilon\,, \\  
&  &
 \gamma_0\leq \frac{n}{2 } -\sqrt{\frac{n^2}{4 }-m^2}
\quad  if \quad 0< m < \frac{n}{2 }, 
\quad while \quad \\
&  &
0 \leq \gamma_0 \leq \frac{n}{2 }  \quad  if \quad \frac{n}{2 } \leq  m \,, 
\end{eqnarray*}
the integral equation (\ref{4.5}) has a unique solution \, $ \psi   (x ,t) \in X({2\varepsilon ,s,\gamma})  $ with 
\[
\cases{ 
\dsp 0 < \gamma <    \frac{1}{\alpha +1}\gamma_0 \quad \mbox{if} \quad 0< m< \frac{n}{2 }\,,\cr
\dsp  
\gamma    \leq  \min \left\{ \gamma _0 \,,\,  \frac{n}{2(\alpha +1)} \right\} 
\quad  if \quad \frac{n}{2 } \leq   m \,.  }
\]
Thus, for the solution $\psi  $  the following estimate holds:
\begin{eqnarray*}  
\sup_{t \in [0,\infty)}  e^{\gamma t}  \|\psi   (\cdot ,t) \|_{H_{(s)} ({\mathbb R}^n)}  < 2\varepsilon \,.
\end{eqnarray*}
\end{theorem}
\medskip

\noindent
{\it Proof.} 
Consider the mapping
\begin{eqnarray} 
\label{5.6}
 S[\psi  ] (x,t)
& := &
\psi _0(x,t) + 
G[ F(\psi   )] (x,t)    \,. 
\end{eqnarray} 
We are going to prove that $S$ maps $X({R,s,\gamma})$ into itself and is a contraction provided that   $\varepsilon  $ and $R$ are sufficiently small. 
\medskip

\noindent
{\bf The case of small physical mass, $  m\in (0, \frac{n}{2 })$.} In this case the operator ${\mathcal K}$ is given by (\ref{operKsm}) and $M=\sqrt{\frac{n^2}{4}-m^2}>0$. 
Corollary~\ref{C3.4}  with $s>n/2$, \,$\gamma =\frac{1}{\alpha +1}(\frac{n}{2}-M - \delta  )>0$\, and \,$\delta   >0 $\,  implies 
\begin{eqnarray*} 
\|  S[\psi  ]  (\cdot  ,t) \|_{ H_{(s)}  ({\mathbb R}^n)  } 
&  \leq   &
 \|\psi  _0(\cdot ,t) \|_{H_{(s)} ({\mathbb R}^n)  }  +     
    \|   G[ F(\cdot , \psi   )](\cdot ,t)  \|_{H_{(s)} ({\mathbb R}^n)  }   \\ 
&  \leq   &
 \|\psi  _0(\cdot ,t) \|_{H_{(s)} ({\mathbb R}^n)  }  + C_M   e^{-(\frac{n}{2}-M) t}    
\int_{ 0}^{t} e^{(\frac{n}{2}-M) b}      \|  F( \cdot ,\psi   )(\cdot ,b)  \|_{H_{(s)} ({\mathbb R}^n)  }  \,db \\
&  \leq   &
   \|\psi  _0(\cdot ,t) \|_{ H_{(s)}({\mathbb R}^n)  }  +  C_M   e^{-\gamma (\alpha +1)  t- \delta  t }    
\int_{ 0}^{t} e^{\gamma (\alpha +1) b+\delta  b}      \|   F( \cdot ,\psi  ) (\cdot ,b) \|_{H_{(s)}({\mathbb R}^n)  }  \,db \,.
\end{eqnarray*}
Taking into account the Condition (${\mathcal L}$)  
we arrive at 
\begin{eqnarray*} 
\|   S[\psi ] (\cdot  ,t) \|_{ H_{(s)}  ({\mathbb R}^n)  } 
&  \leq   &
  \|\psi  _0(\cdot ,t) \|_{ H_{(s)}({\mathbb R}^n)  }  +  C_M   e^{-\gamma (\alpha +1) t- \delta  t }    
\int_{ 0}^{t} e^{\gamma (\alpha +1) b+\delta  b}     \| \psi  (\cdot ,b ) \| _{H_{(s)}({\mathbb R}^n)}   ^{\alpha  +1}  \,db  \\
&  \leq   &
  \|\psi  _0(\cdot ,t) \|_{ H_{(s)} ({\mathbb R}^n)  }  +  C_M   e^{-\gamma (\alpha +1) t- \delta  t }    
\int_{ 0}^{t} e^{ \delta  b} \left( e^{ \gamma  b }     \| \psi  (\cdot ,b ) \| _{H_{(s)}({\mathbb R}^n)}   \right)^{\alpha  +1}  \,db  \,  .  
\end{eqnarray*}
Then 
\begin{eqnarray*} 
&  &
e^{\gamma t}\|   S[\psi  ] (x ,t) \|_{ H_{(s)}({\mathbb R}^n)  } \quad \leq \quad
e^{\gamma (\alpha +1) t}\|   S[\psi  ] (x ,t) \|_{ H_{(s)}({\mathbb R}^n)  } \\
&  \leq   &
 e^{\gamma (\alpha +1) t}\|\psi  _0(\cdot ,t) \|_{H_{(s)}({\mathbb R}^n)  }  +  C_M       
 \left( \sup_{\tau  \in [0,\infty)} e^{\gamma  \tau  }     \| \psi   (\cdot ,\tau  ) \| _{H_{(s)}({\mathbb R}^n)}   \right)^{\alpha  +1} e^{- \delta  t }\int_{ 0}^{t} e^{ \delta  b} \,db \\
&  \leq   &
 e^{\gamma_0 t}\|\psi  _0(\cdot ,t) \|_{H_{(s)}({\mathbb R}^n)  }  +  C_M  \delta ^{-1}     
 \left( \sup_{\tau  \in [0,\infty)} e^{\gamma   \tau  }     \| \psi   (\cdot ,\tau  ) \| _{H_{(s)}({\mathbb R}^n)}   \right)^{\alpha  +1} ,  
\end{eqnarray*}
and
\begin{eqnarray*} 
&  &
\sup_{t  \in [0,\infty)}e^{\gamma t}\|   S[\psi  ] (x ,t) \|_{ H_{(s)}({\mathbb R}^n)  } \\
&  \leq   &
 \sup_{t  \in [0,\infty)}e^{\gamma_0 t}\|\psi  _0(\cdot ,t) \|_{H_{(s)}({\mathbb R}^n)  }  +  C_M  \delta ^{-1}     
 \left( \sup_{t  \in [0,\infty)} e^{\gamma   t  }     \| \psi   (\cdot ,t  ) \| _{H_{(s)}({\mathbb R}^n)}   \right)^{\alpha  +1} .  \nonumber
\end{eqnarray*}
In particular, if $\gamma_0 = \frac{n}{2}-M >0$, then, with $\delta >0 $ such that $ \gamma (\alpha +1)= \frac{n}{2}-M - \delta < \gamma_0$, we have 
\begin{eqnarray*}  
&  &
\sup_{t \in [0,\infty)}  e^{\gamma t} \|S[\psi  ](\cdot,t) \|_{H_{(s)}({\mathbb R}^n)} \\ 
& \leq  &
\sup_{t \in [0,\infty)}  e^{( \frac{n}{2}-M)  t} \|\psi  _0(\cdot ,t) \|_{H_{(s)}({\mathbb R}^n)} + C\left(\sup_{t \in [0,\infty)}  e^{\gamma  t}   
 \| \psi   (\cdot ,t ) \| _{H_{(s)}({\mathbb R}^n)}  \right) ^{\alpha  +1} \,.       
\end{eqnarray*}
 Thus, the last inequality proves that the operator $S$ maps $X({R,s,\gamma})$ into itself if $\varepsilon  $ and $R$ are sufficiently small, namely, if
$\varepsilon  +C R^{\alpha +1} < R $. 
\medskip

It remains to prove that $S$ is a contraction mapping.
As a matter of fact, we just need to apply estimate (\ref{calM})     and get the contraction property from 
\[
 e^{\gamma t} \|S[\psi  ](\cdot,t) -  S[\widetilde{\psi }   ](\cdot,t) \|_{H_{(s)}({\mathbb R}^n) }
 \leq  
CR(t) ^{\alpha } d(\psi  ,\widetilde{\psi } )\,,  
\]
where $\displaystyle   R(t):= \max\{ \sup_{0\leq \tau \leq t } e^{\gamma \tau  }  \| \psi  (\cdot ,\tau ) \| _{H_{(s)}({\mathbb R}^n) } , \sup_{0\leq \tau \leq t } e^{\gamma \tau  }  \| \widetilde{\psi }   (\cdot ,\tau ) \| _{H_{(s)}({\mathbb R}^n)  }\} \leq R$. 
Indeed, due to Condition~(${\mathcal L}$),  we have
\begin{eqnarray*}
&  &
\| S[\psi  ](\cdot , t) -  S[\widetilde{\psi }  ](\cdot , t) \|_{H_{(s)}({\mathbb R}^n) }  
 = 
\| G[ \,
( F(\cdot , \psi   ) - F( \cdot ,\widetilde{\psi }  ) )  ](\cdot, t)
\|_{H_{(s)}({\mathbb R}^n) }  \\
& \le  &
C_M     e^{-(\frac{n}{2}-M) t}  
\int_{ 0}^{t}  e^{(\frac{n}{2}-M) b} \|  ( F( \cdot ,\psi   ) - F(\cdot , \widetilde{\psi } ) )  (\cdot ,b)) \| _{H_{(s)}({\mathbb R}^n) } \,db \\
& \le  &
C_M     e^{-\gamma (\alpha +1) t- \delta  t }   
\int_{ 0}^{t}  e^{\gamma (\alpha +1) b+\delta  b} \|  ( F( \cdot ,\psi   ) - F( \cdot ,\widetilde{\psi } ) )  (\cdot ,b)) \| _{H_{(s)}({\mathbb R}^n) } \,db    \\
& \le  &
C_M     e^{-\gamma (\alpha +1) t- \delta  t }   
\int_{ 0}^{t}  e^{\gamma (\alpha +1) b+\delta  b} \|\psi  (\cdot , b) -\widetilde{\psi }(\cdot , b) \|_{H_{(s)}({\mathbb R}^n) } \\
&  &
\hspace{3.5cm}\times \Big( \| \psi   (\cdot ,  b) \|_{H_{(s)}({\mathbb R}^n) } ^\alpha  
+ \| \widetilde{\psi }(\cdot ,  b)  \|_{H_{(s)}({\mathbb R}^n) }^\alpha 
\Big) \,db    \,.
\end{eqnarray*}
Thus,  taking into account  the last estimate and the definition of the metric $ d(\psi  ,\widetilde{\psi }) $, we obtain
\begin{eqnarray*}
&  &
e^{\gamma (\alpha +1) t }  \| S[\psi  ](\cdot , t) -  S[\widetilde{\psi } ](\cdot , t) \|_{H_{(s)}({\mathbb R}^n) }   \\ 
& \le  &
C_M     e^{- \delta  t }   
\int_{ 0}^{t}  e^{\gamma (\alpha +1) b+\delta  b} \|\psi  (\cdot , b) -\widetilde{\psi } (\cdot , b) \|_{H_{(s)}({\mathbb R}^n) } 
\Big( \| \psi   (\cdot ,  b) \|_{H_{(s)}({\mathbb R}^n) } ^\alpha  
+ \| \widetilde{\psi }(\cdot ,  b)  \|_{H_{(s)}({\mathbb R}^n) }^\alpha 
\Big) \,db   \\ 
& \le  &
C_M     e^{- \delta  t }   
\int_{ 0}^{t}  e^{\delta  b}\Big( \max_{0 \le \tau  \leq b } e^{\gamma \tau }\|\psi  (\cdot , \tau ) -\widetilde{\psi } (\cdot , \tau ) \|_{H_{(s)}({\mathbb R}^n) }  \Big)\\
&  &
\times 
\Big(  \Big(  \max_{0 \le \tau  \leq b }e^{\gamma  \tau }\| \psi   (\cdot ,  \tau ) \|_{H_{(s)}({\mathbb R}^n) }  \Big)^\alpha  
+  \Big(  \max_{0 \le \tau  \leq b }e^{\gamma  \tau }\| \widetilde{\psi }  (\cdot ,  \tau ) \|_{H_{(s)}({\mathbb R}^n) }  \Big)^\alpha   
\Big) \,db   \\
& \le  &
C_{M,\alpha}   d(\psi  ,\widetilde{\psi }) R(t)^\alpha      
 e^{- \delta  t }   
\int_{ 0}^{t}  e^{\delta  b}   \,db   \\
& \le &
 C_{M,\alpha} \delta^{-1}  d(\psi  ,\widetilde{\psi } ) R(t)^\alpha   \,.
\end{eqnarray*}
Consequently, 
\begin{eqnarray*}
e^{\gamma  t }  \| S[\psi  ](\cdot , t) -  S[\widetilde{\psi } ](\cdot , t) \|_{H_{(s)}({\mathbb R}^n) }  
& \le  &
 C_{M,\alpha} \delta^{-1} R(t)^\alpha   d(\psi  ,\widetilde{\psi } )  \,.
\end{eqnarray*}
Set $R:=\sup_{t \in [0,\infty)} R(t) $. Then we choose $\varepsilon $ and $ R$ such that $C_{M,\alpha} \delta^{-1}   R ^\alpha <1 $.  Banach's fixed point theorem completes the proof for the case of  small 
physical mass.
\medskip

\noindent
{\bf The case of the  large physical mass $m  \geq n /2 $.} In this case the operator ${\mathcal K}$ is given by (\ref{operK}). We set $\gamma \leq \min \{\gamma _0\,,\,\frac{n}{2(\alpha +1)} \}$ in the definition of metric of the space $X({R,s,\gamma})$. Then, due to (\ref{2.13}),  we have
\begin{eqnarray*} 
\|    S[\psi  ] (\cdot  ,t) \|_{ H_{(s)}  ({\mathbb R}^n)  } 
&  \leq   &
 \|\psi  _0(\cdot ,t) \|_{H_{(s)} ({\mathbb R}^n)  }  +     
    \|   G[ F( \cdot ,\psi   )](\cdot ,t)  \|_{H_{(s)} ({\mathbb R}^n)  }   \\ 
&  \leq   &
 \|\psi  _0(\cdot ,t) \|_{H_{(s)} ({\mathbb R}^n)  }  + C_M   e^{-\frac{n}{2} t}    
\int_{ 0}^{t} e^{\frac{n}{2} b}   (1+ t -b )^{1- \sgn M}   \| F( \psi   )(\cdot ,b)  \|_{H_{(s)} ({\mathbb R}^n)  }  \,db \\
&  \leq   &
 \|\psi  _0(\cdot ,t) \|_{H_{(s)} ({\mathbb R}^n)  }  + C_{M,\alpha}   e^{-\frac{n}{2} t}    
\int_{ 0}^{t} e^{\frac{n}{2} b}    (1+ t -b )^{1- \sgn M}  \| \psi   (\cdot ,b ) \| _{H_{(s)}({\mathbb R}^n)}   ^{\alpha  +1}   \,db \,.   
\end{eqnarray*}
Hence
\begin{eqnarray*} 
&  &
e^{\gamma t}\|    S[\psi  ] (\cdot  ,t) \|_{ H_{(s)}({\mathbb R}^n)  } \\
&  \leq   &
e^{\gamma t}\|\psi  _0(\cdot ,t) \|_{H_{(s)}({\mathbb R}^n)  }  \\
&  &
+  C_{M,\alpha}       
 \left( \sup_{\tau  \in [0,\infty)} e^{\gamma \tau } \| \psi   (\cdot ,\tau  ) \| _{H_{(s)}({\mathbb R}^n)}   \right)^{\alpha  +1} e^{-(\frac{n}{2}-\gamma )t}    
\int_{ 0}^{t} e^{(\frac{n}{2}-\gamma (\alpha +1))b}   (1+ t -b )^{1- \sgn M}  \,db   \\
&  \leq   &
e^{\gamma_0 t}\|\psi  _0(\cdot ,t) \|_{H_{(s)}({\mathbb R}^n)  }   
+  C_{M,\alpha}       
 \left( \sup_{\tau  \in [0,\infty)} e^{\gamma \tau } \| \psi   (\cdot ,\tau  ) \| _{H_{(s)}({\mathbb R}^n)}   \right)^{\alpha  +1} e^{-\gamma \alpha t}    .  
\end{eqnarray*}
Then we choose $\varepsilon $ and $ R$ such that $\varepsilon +  C_{M,\alpha}   R ^{\alpha+1} <R $.

To prove that $S$ is a contraction mapping, we just need to apply estimate (\ref{calM})     and get the contraction property from 
\[
e^{\gamma t} \|S[\psi  ](\cdot,t) -  S[\widetilde{\psi } ](\cdot,t) \|_{H_{(s)}({\mathbb R}^n) }
 \leq  
CR(t) ^{\alpha } d(\psi  ,\widetilde{\psi } )\,,  
\]
where $\displaystyle R(t):= \max\{ \sup_{0\leq \tau \leq t }   e^{\gamma \tau }\| \psi   (\cdot ,\tau ) \| _{H_{(s)}({\mathbb R}^n) }
, \sup_{0\leq \tau \leq t }  e^{\gamma \tau } \| \widetilde{\psi }(\cdot ,\tau ) \| _{H_{(s)}({\mathbb R}^n) } \}\leq R$. 
Indeed,   we have
\begin{eqnarray*}
&  &
e^{\gamma t}\| S[\psi  ](\cdot , t) -  S[\widetilde{\psi } ](\cdot , t) \|_{H_{(s)}({\mathbb R}^n) }  
 = 
e^{\gamma t}\| G[ \,
( F( \cdot ,\psi   ) - F( \cdot ,\widetilde{\psi } ) )  ](\cdot, t)
\|_{H_{(s)}({\mathbb R}^n) }  \\
& \le  &
C_M     e^{-(\frac{n}{2}-\gamma ) t}  
\int_{ 0}^{t}  e^{\frac{n}{2} b}  (1+ t -b )^{1- \sgn M}\|  ( F(\cdot , \psi   ) - F(\cdot , \widetilde{\psi } ) )  (\cdot ,b)) \| _{H_{(s)}({\mathbb R}^n) } \,db \\
& \le  &
C_{M,\alpha }    e^{-(\frac{n}{2}-\gamma ) t}    
\int_{ 0}^{t}  e^{\frac{n}{2} b} (1+ t -b )^{1- \sgn M}\|\psi  (\cdot , b) -\widetilde{\psi }(\cdot , b) \|_{H_{(s)}({\mathbb R}^n) } \\
&  &
\times 
\Big( \| \psi   (\cdot ,  b) \|_{H_{(s)}({\mathbb R}^n) } ^\alpha  
+ \| \widetilde{\psi }  (\cdot ,  b)  \|_{H_{(s)}({\mathbb R}^n) }^\alpha 
\Big) \,db    \,.
\end{eqnarray*}
Thus,  taking into account  the last estimate  and a definition of the metric, we obtain
\begin{eqnarray*}
&  &
 e^{\gamma t}\| S[\psi  ](\cdot , t) -  S[\Psi ](\cdot , t) \|_{H_{(s)}({\mathbb R}^n) }   \\ 
& \le  &
C_{M,\alpha}      e^{- (\frac{n}{2}-\gamma )  t }   
\int_{ 0}^{t}  e^{(\frac{n}{2}-\gamma (\alpha +1) ) b} (1+ t -b )^{1- \sgn M}e^{\gamma   b}\|\psi  (\cdot , b) -\widetilde{\psi } (\cdot , b) \|_{H_{(s)}({\mathbb R}^n) } \\
&  &
\times \Big( (e^{\gamma   b}\| \psi   (\cdot ,  b) \|_{H_{(s)}({\mathbb R}^n) } )^\alpha  
+ (e^{\gamma   b}\| \widetilde{\psi } (\cdot ,  b)  \|_{H_{(s)}({\mathbb R}^n) })^\alpha 
\Big) \,db   \\ 
& \le  &
C_{M,\alpha}   d(\psi  ,\widetilde{\psi }) R(t)^\alpha      
 e^{- (\frac{n}{2}-\gamma )  t }   
\int_{ 0}^{t}   e^{(\frac{n}{2}-\gamma (\alpha +1) ) b}(1+ t -b )^{1- \sgn M}  \,db    \,,
\end{eqnarray*}
and, consequently, 
\begin{eqnarray*}
 e^{\gamma t}\| S[\psi  ](\cdot , t) -  S[\widetilde{\psi }](\cdot , t) \|_{H_{(s)}({\mathbb R}^n) }  
& \le &
 C_{M,\alpha}    d(\psi  ,\widetilde{\psi } ) R(t)^\alpha  e^{-  \gamma  \alpha   t }  \,.
\end{eqnarray*}
Set $R:=\sup_{t \in [0,\infty)} R(t) $. Then we choose $\varepsilon $ and $ R$ such that $ C_{M,\alpha}    R ^\alpha  <1 $. 
Banach's fixed point theorem completes the proof of theorem.
\hfill $\square$

\subsection{Proof of Theorem~\ref{T0.1}} 
\label{SS5.2}
 
{\bf The case of the small physical mass inside of Higuchi bound, $m  <\sqrt{n^2-1}/2  $.} In this case the operator ${\mathcal K}$ is given by (\ref{operKsm}) and $M=\sqrt{\frac{n^2}{4}-m^2}$. 
Then for the function   $\psi _0=  \psi _0(x,t)$, that is, for the solution of the Cauchy problem (\ref{1.10})  and for $s>\frac{n}{2}$, $p=p'=2$, $n \geq 2$, according to Theorem~\ref{T3.1}
we have the estimate 
\begin{equation}
\label{5.7} 
\| \psi _0 (\cdot ,t) \|_{ { H}_{(s)} ({\mathbb R}^n)}
  \leq   
C_{M,n,p,q,s}  e^{( M -\frac{n}{2})t}\Big\{   \| \psi _0   \|_{ { H}_{(s)} ({\mathbb R}^n)}
+  \|\psi _1  \|_{ { H}_{(s)} ({\mathbb R}^n)} 
 \Big\} \,.
\end{equation}
For every initial functions $\psi  _0 (x)$ and $\psi  _1 (x)$ 
the function $\psi =\psi _0 (x,t)$ belongs to the space  $X(R,s,\gamma ) $, where the operator $S $ (\ref{5.6}) is a contraction.  
In the case of $n=3$ that means $m^2 < 2 $.
Theorem~\ref{TIE} completes the proof of the existence of the global solution.
\medskip

\noindent
{\bf The case of the critical mass}, $ m=\sqrt{n^2-1}/2$. We apply Lemma~\ref{L1} that shows that the estimate (\ref{5.7})  holds with $ M=1/2$.  

\medskip

\noindent
{\bf The case of the small physical mass outside of Higuchi bound, $\sqrt{n^2-1}/2<m< n/2 $.} In this case  $\psi  _0 (x)=0 $ and again, 
the function $\psi =\psi _0 (x,t)$ belongs to the space  $X(R,s,\gamma ) $, where the operator $S $ (\ref{5.6}) is a contraction.
\medskip

\noindent
{\bf The case of the large physical mass $m \geq  n /2 $.} In this case the operator ${\mathcal K}$ is given by (\ref{operK}).  Then for the function   $\psi  _0 $, that is for the solution of the Cauchy problem (\ref{1.10})  and for $s>\frac{n}{2}$, $p=p'=2$, $n \geq 2$ we apply the estimate (\ref{2.4}), 
\begin{eqnarray*}
e^{\gamma_0 t}\|  \psi _0 (\cdot ,t) \|_{ { H}_{(s)} ({\mathbb R}^n)  } 
&  \leq   &
C_M e^{\gamma_0 t-\frac{n}{2}t}(1+ t )^{1- \sgn M} \Big\{ e ^{\frac{t}{2}}  \| \psi _0  \|_{ { H}_{(s)} ({\mathbb R}^n)  }
+  \|\psi _1  \|_{ { H}_{(s)} ({\mathbb R}^n)  } 
 \Big\} \\
& \leq  &
C_M e^{(\gamma_0  -\frac{n-1}{2})t}(1+ t )^{1- \sgn M} \Big\{  \| \psi _0  \|_{ { H}_{(s)} ({\mathbb R}^n)  }
+  \|\psi _1  \|_{ { H}_{(s)} ({\mathbb R}^n)  } 
 \Big\}\,.
\end{eqnarray*}
We choose $0\leq \gamma _0 <  \frac{n-1}{2}$ if  $m=n/2$ and $0\leq \gamma _0 \leq   \frac{n-1}{2}$ if $m>n/2$.
Thus,  $ \psi _0 \in X(R,s,\gamma _0)$. Theorem~\ref{TIE}  implies   existence of a solution  $ \psi   \in X(R,s,\gamma )$  with  
$\gamma    \leq  \min \left\{ \gamma _0 \,,\,  \frac{n}{2(\alpha +1)} \right\} $
of the integral equation (\ref{5.5})  provided that $R$ is 
sufficiently small.  Theorem~\ref{T0.1} is proved. \hfill $\square$
\medskip

\subsection{Proof of Theorem~\ref{T0.2}} 
\label{SS5.3}

\noindent
 First we consider the case of small  { $m<n/2$}.  According to Corollary~\ref{C3.3} the solution $\psi _0 $ of the linear problem (\ref{2.5}) 
 with ``plus'' 
 satisfies the following  estimate 
\begin{eqnarray*} 
\| \psi_0 (x ,t) \|_{H_{(s)} ({\mathbb R}^n)} 
&  \leq   &
 C_M   e^{-(\frac{n}{2}-M)t}    
\int_{ 0}^{t} e^{(\frac{n}{2}-M)b}  \| f(x ,b)  \|_{H_{(s)} ({\mathbb R}^n)} \,db \\
&  \leq   &
 C_M   e^{-(\frac{n}{2}-M)t}    
\int_{ 0}^{t} e^{(\frac{n}{2}-M)b}  e^{-\gamma _{rhs}b } e^{\gamma _{rhs} b} \| f(x ,b)  \|_{H_{(s)} ({\mathbb R}^n)} \,db \\
&  \leq   &
 C_M     
\left( \sup_{\tau  \in [0,t)} e^{\gamma _{rhs} \tau }\| f (x ,\tau ) \|_{H_{(s)} ({\mathbb R}^n)} \right)  e^{-(\frac{n}{2}-M)t} \int_{ 0}^{t} e^{(\frac{n}{2}-M-\gamma _{rhs})b}   \,db\,.  
\end{eqnarray*}
Consider three cases: $\gamma _{rhs} < n/2-M $, $\gamma _{rhs} = n/2-M$, and $\gamma _{rhs} >n/2-M$. \\
In the first case of   { $\gamma _{rhs} < n/2-M $} we have 
\begin{eqnarray*} 
\| \psi_0 (x ,t) \|_{H_{(s)} ({\mathbb R}^n)} 
&  \leq   &
 C_M     
\left( \sup_{\tau  \in [0,\infty)} e^{\gamma _{rhs} \tau }\| f (x ,\tau ) \|_{H_{(s)} ({\mathbb R}^n)} \right)  e^{-\gamma _{rhs} t}  \,.
\end{eqnarray*}
If  {  $\gamma _{rhs} = n/2-M$}, then
\begin{eqnarray*} 
e^{(\frac{n}{2}-M)t}  \| \psi_0 (x ,t) \|_{H_{(s)} ({\mathbb R}^n)} 
&  \leq   &
 C_M t    
\left( \sup_{\tau  \in [0,\infty)} e^{\gamma _{rhs} \tau }\| f (x ,\tau ) \|_{H_{(s)} ({\mathbb R}^n)} \right)  \,.  
\end{eqnarray*}
Thus, according to Theorem~\ref{TIE} for sufficiently small $ \varepsilon $, and 
\begin{eqnarray*}  
 \sup_{\tau  \in [0,\infty)} e^{\gamma _{rhs} \tau }\| f (x ,\tau ) \|_{H_{(s)} ({\mathbb R}^n)}  \leq \varepsilon    
\end{eqnarray*}
 the problem (\ref{NWEf}), (\ref{ICPHIf}) has a global solution and 
\begin{eqnarray*}  
  \sup_{\tau  \in [0,\infty)} e^{\gamma   \tau }\| \psi  (x ,\tau ) \|_{H_{(s)} ({\mathbb R}^n)}  \leq 2\varepsilon \,,   
\end{eqnarray*} 
where $ \gamma < \gamma _{rhs}/(\alpha +1)$.

If   { $\gamma _{rhs} > n/2-M$}, then 
\begin{eqnarray*} 
e^{(\frac{n}{2}-M)t} \| \psi_0 (x ,t) \|_{H_{(s)} ({\mathbb R}^n)} 
&  \leq   &
 C_M     t
\left( \sup_{\tau  \in [0,\infty)} e^{\gamma _{rhs} \tau }\| f (x ,\tau ) \|_{H_{(s)} ({\mathbb R}^n)} \right)  ,  
\end{eqnarray*}
and 
\begin{eqnarray*}  
  \sup_{\tau  \in [0,\infty)} e^{\gamma   \tau }\| \psi  (x ,\tau ) \|_{H_{(s)} ({\mathbb R}^n)}  \leq 2\varepsilon \,,   
\end{eqnarray*} 
where $ \gamma <(\frac{n}{2}-M)/(\alpha +1)$. 

For the {large mass  $m\geq n/2 $}, due to (\ref{2.13}),  we have
\begin{eqnarray*} 
&  &
\| \psi_0 (x ,t) \|_{H_{(s)} ({\mathbb R}^n)} \\
&  \leq   &
 C_M   e^{- \frac{n}{2} t}    
\int_{ 0}^{t} e^{ \frac{n}{2} b}  \| f(x ,b)  \|_{H_{(s)} ({\mathbb R}^n)} (1+t-b)^{1-\sgn M}\,db \\
&  \leq   &
 C_M  \left( \sup_{\tau  \in [0,\infty)} e^{\gamma _{rhs} \tau }\| f (x ,\tau ) \|_{H_{(s)} ({\mathbb R}^n)} \right)  e^{- \frac{n}{2} t}    
\int_{ 0}^{t} e^{(\frac{n}{2}-\gamma _{rhs})b} (1+t-b)^{1-\sgn M}\,db \,.
\end{eqnarray*}
First let { $m=n/2$}, then
\begin{eqnarray*}  
\| \psi_0 (x ,t) \|_{H_{(s)} ({\mathbb R}^n)}  
&  \leq   &
 C_M  \left( \sup_{\tau  \in [0,\infty)} e^{\gamma _{rhs} \tau }\| f (x ,\tau ) \|_{H_{(s)} ({\mathbb R}^n)} \right)  e^{- \frac{n}{2} t}    
\int_{ 0}^{t} e^{(\frac{n}{2}-\gamma _{rhs})b} (1+t-b)\,db \,.
\end{eqnarray*}
If,  in addition, { $\frac{n}{2}>\gamma _{rhs} $}, then 
\begin{eqnarray*} 
\| \psi_0 (x ,t) \|_{H_{(s)} ({\mathbb R}^n)}  
&  \leq   &
 C_M  \left( \sup_{\tau  \in [0,\infty)} e^{\gamma _{rhs} \tau }\| f (x ,\tau ) \|_{H_{(s)} ({\mathbb R}^n)} \right)      
 e^{ -\gamma _{rhs}t} 
\end{eqnarray*}
and we choose $\gamma _0=\gamma _{rhs} $  and apply Theorem~\ref{TIE} which implies existence of a global solution $\psi  \in X(R,s,\gamma ) $, where
$\gamma  \leq \min\{\gamma _{rhs}, \frac{n}{2(\alpha +1)}\} $. \\
If additionally { $\frac{n}{2}=\gamma _{rhs} $}, then 
\begin{eqnarray*} 
\| \psi_0 (x ,t) \|_{H_{(s)} ({\mathbb R}^n)}  
&  \leq   &
 C_M  \left( \sup_{\tau  \in [0,\infty)} e^{\gamma _{rhs} \tau }\| f (x ,\tau ) \|_{H_{(s)} ({\mathbb R}^n)} \right)    
 e^{ -\gamma _{rhs}t} (1+t)^2\,.  
\end{eqnarray*}
We choose $\gamma _0<\gamma _{rhs} $  and apply Theorem~\ref{TIE} which implies existence of a global solution $\psi  \in X(R,s,\gamma ) $, where
$\gamma  \leq \min\{\gamma _{0}, \frac{n}{2(\alpha +1)}\} $. \\
If additionally{  $\frac{n}{2}<\gamma _{rhs} $}, then
\begin{eqnarray*} 
\| \psi_0 (x ,t) \|_{H_{(s)} ({\mathbb R}^n)} 
&  \leq   &
 C_M  \left( \sup_{\tau  \in [0,\infty)} e^{\gamma _{rhs} \tau }\| f (x ,\tau ) \|_{H_{(s)} ({\mathbb R}^n)} \right)  e^{- \frac{n}{2} t}    
\int_{ 0}^{t} e^{(\frac{n}{2}-\gamma _{rhs})b} (1+t-b)\,db \\
&  \leq   &
 C_M  \left( \sup_{\tau  \in [0,\infty)} e^{\gamma _{rhs} \tau }\| f (x ,\tau ) \|_{H_{(s)} ({\mathbb R}^n)} \right)  e^{- \frac{n}{2} t} (1+t)
\end{eqnarray*}
and we choose $\gamma _0<n/2 $ and $\gamma \leq \frac{n}{2(\alpha +1)} $.
\smallskip

\noindent
Next we consider the case of {$m>n/2 $}, then
\begin{eqnarray*}  
\| \psi_0 (x ,t) \|_{H_{(s)} ({\mathbb R}^n)}  
&  \leq   &
 C_M  \left( \sup_{\tau  \in [0,\infty)} e^{\gamma _{rhs} \tau }\| f (x ,\tau ) \|_{H_{(s)} ({\mathbb R}^n)} \right)  e^{- \frac{n}{2} t}    
\int_{ 0}^{t} e^{(\frac{n}{2}-\gamma _{rhs})b}  \,db \,.
\end{eqnarray*}
If additionally{  $\frac{n}{2}>\gamma _{rhs} $}, then
\begin{eqnarray*}  
\| \psi_0 (x ,t) \|_{H_{(s)} ({\mathbb R}^n)}  
&  \leq   &
 C_M  \left( \sup_{\tau  \in [0,\infty)} e^{\gamma _{rhs} \tau }\| f (x ,\tau ) \|_{H_{(s)} ({\mathbb R}^n)} \right)  e^{- \gamma _{rhs} t} 
\end{eqnarray*}
and we choose $\gamma _0=\gamma _{rhs}  $ and $\gamma \leq \min \{ \gamma _{rhs}, \frac{n}{2(\alpha +1)}\} $.\\
If additionally { $\frac{n}{2}=\gamma _{rhs} $}, then 
\begin{eqnarray*}  
\| \psi_0 (x ,t) \|_{H_{(s)} ({\mathbb R}^n)}  
&  \leq   &
 C_M  \left( \sup_{\tau  \in [0,\infty)} e^{\gamma _{rhs} \tau }\| f (x ,\tau ) \|_{H_{(s)} ({\mathbb R}^n)} \right)  e^{- \frac{n}{2} t} t
\end{eqnarray*}
and we choose $\gamma _0<n/2 $ and $\gamma \leq \min\{ \frac{n}{2}\,,\, \frac{n}{2(\alpha +1)} \}=  \frac{n}{2(\alpha +1)}$.\\
If additionally{  $\frac{n}{2}<\gamma _{rhs} $}, then
\begin{eqnarray*}  
\| \psi_0 (x ,t) \|_{H_{(s)} ({\mathbb R}^n)}  
&  \leq   &
 C_M  \left( \sup_{\tau  \in [0,\infty)} e^{\gamma _{rhs} \tau }\| f (x ,\tau ) \|_{H_{(s)} ({\mathbb R}^n)} \right)  e^{- \frac{n}{2} t}   
\end{eqnarray*}
and we choose $\gamma _0=n/2 $ and    $\gamma \leq \min\{ \frac{n}{2}\,,\, \frac{n}{2(\alpha +1)} \}=  \frac{n}{2(\alpha +1)}$.
 Theorem~\ref{T0.2} is proved. \hfill $\square$

\end{document}